\documentclass{amsart} 
\usepackage{amssymb, amsbsy, amsthm, amsmath,amstext, amsopn} 
\usepackage[all]{xy} 
\usepackage{amsfonts}
\usepackage{amscd} 
\newtheorem{thm}{Theorem} [section]
\newtheorem{lemma}[thm]{Lemma}

\newtheorem{corollary}[thm]{Corollary}
\newtheorem{prop}[thm]{Proposition}

\theoremstyle{definition}

\newtheorem{defn}[thm]{Definition}

\newtheorem{example}[thm]{Example}

\theoremstyle{remark}

\newtheorem{remark}[thm]{Remark}

\begin{document}

\numberwithin{equation}{section}

\newcommand{\hs}{\mbox{\hspace{.4em}}}
\newcommand{\ds}{\displaystyle}
\newcommand{\bd}{\begin{displaymath}}
\newcommand{\ed}{\end{displaymath}}
\newcommand{\bcd}{\begin{CD}}
\newcommand{\ecd}{\end{CD}}

\newcommand{\on}{\operatorname}
\newcommand{\proj}{\operatorname{Proj}}
\newcommand{\bproj}{\underline{\operatorname{Proj}}}
\newcommand{\spec}{\operatorname{Spec}}
\newcommand{\bspec}{\underline{\operatorname{Spec}}}
\newcommand{\pline}{{\mathbf P} ^1}
\newcommand{\aline}{{\mathbf A} ^1}
\newcommand{\pplane}{{\mathbf P}^2}
\newcommand{\aplane}{{\mathbf A}^2}
\newcommand{\coker}{{\operatorname{coker}}}
\newcommand{\ldb}{[[}
\newcommand{\rdb}{]]}

\newcommand{\Sym}{\operatorname{Sym}^{\bullet}}
\newcommand{\Symp}{\operatorname{Sym}}
\newcommand{\Pic}{\operatorname{Pic}}
\newcommand{\AAut}{\operatorname{Aut}}
\newcommand{\PAut}{\operatorname{PAut}}

\newcommand{\too}{\twoheadrightarrow}
\newcommand{\C}{{\mathbf C}}
\newcommand{\Z}{{\mathbf Z}}
\newcommand{\R}{{\mathbf R}}
\newcommand{\Cx}{{\mathbf C}^{\times}}
\newcommand{\Cbar}{\overline{\C}}
\newcommand{\Cxbar}{\overline{\Cx}}
\newcommand{\cA}{{\mathcal A}}
\newcommand{\cS}{{\mathcal S}}
\newcommand{\cV}{{\mathcal V}}
\newcommand{\cM}{{\mathcal M}}
\newcommand{\bA}{{\mathbf A}}
\newcommand{\cB}{{\mathcal B}}
\newcommand{\cC}{{\mathcal C}}
\newcommand{\cD}{{\mathcal D}}
\newcommand{\D}{{\mathcal D}}
\newcommand{\cs}{{\mathbf C} ^*}
\newcommand{\boldc}{{\mathbf C}}
\newcommand{\cE}{{\mathcal E}}
\newcommand{\cF}{{\mathcal F}}
\newcommand{\bF}{{\mathbf F}}
\newcommand{\cG}{{\mathcal G}}
\newcommand{\G}{{\mathbb G}}
\newcommand{\cH}{{\mathcal H}}
\newcommand{\cJ}{{\mathcal J}}
\newcommand{\cK}{{\mathcal K}}
\newcommand{\cL}{{\mathcal L}}
\newcommand{\baL}{{\overline{\mathcal L}}}
\newcommand{\M}{{\mathcal M}}
\newcommand{\Mf}{{\mathfrak M}}
\newcommand{\bM}{{\mathbf M}}
\newcommand{\bm}{{\mathbf m}}
\newcommand{\cN}{{\mathcal N}}
\newcommand{\theo}{\mathcal{O}}
\newcommand{\cP}{{\mathcal P}}
\newcommand{\cR}{{\mathcal R}}
\newcommand{\Pp}{{\mathbb P}}
\newcommand{\boldp}{{\mathbf P}}
\newcommand{\boldq}{{\mathbf Q}}
\newcommand{\bbL}{{\mathbf L}}
\newcommand{\cQ}{{\mathcal Q}}
\newcommand{\cO}{{\mathcal O}}
\newcommand{\Oo}{{\mathcal O}}
\newcommand{\OX}{{\Oo_X}}
\newcommand{\OY}{{\Oo_Y}}
\newcommand{\otY}{{\underset{\OY}{\ot}}}
\newcommand{\otX}{{\underset{\OX}{\ot}}}
\newcommand{\cU}{{\mathcal U}}\newcommand{\cX}{{\mathcal X}}
\newcommand{\cW}{{\mathcal W}}
\newcommand{\boldz}{{\mathbf Z}}
\newcommand{\qgr}{\operatorname{q-gr}}
\newcommand{\gr}{\operatorname{gr}}
\newcommand{\rk}{\operatorname{rk}}
\newcommand{\coh}{\operatorname{coh}}
\newcommand{\End}{\operatorname{End}}
\newcommand{\uEnd}{\underline{\operatorname{End}}}
\newcommand{\Hom}{\operatorname{Hom}}
\newcommand{\uHom}{\underline{\operatorname{Hom}}}
\newcommand{\uHomY}{\uHom_{\OY}}
\newcommand{\uHomX}{\uHom_{\OX}}
\newcommand{\Ext}{\operatorname{Ext}}
\newcommand{\bExt}{\operatorname{\bf{Ext}}}
\newcommand{\Tor}{\operatorname{Tor}}

\newcommand{\inv}{^{-1}}
\newcommand{\airtilde}{\widetilde{\hspace{.5em}}}
\newcommand{\airhat}{\widehat{\hspace{.5em}}}
\newcommand{\nt}{^{\circ}}
\newcommand{\del}{\partial}

\newcommand{\supp}{\operatorname{supp}}
\newcommand{\GK}{\operatorname{GK-dim}}
\newcommand{\hd}{\operatorname{hd}}
\newcommand{\id}{\operatorname{id}}
\newcommand{\res}{\operatorname{res}}
\newcommand{\lrar}{\leadsto}
\newcommand{\im}{\operatorname{Im}}
\newcommand{\HH}{\operatorname{H}}
\newcommand{\TF}{\operatorname{TF}}
\newcommand{\Bun}{\operatorname{Bun}}
\newcommand{\BunD}{\operatorname{Bun}_{\D^{\bullet}}}
\newcommand{\PicD}{\operatorname{Pic}_{\D}}
\newcommand{\Hilb}{\operatorname{Hilb}}
\newcommand{\Fact}{\operatorname{Fact}}
\newcommand{\CM}{\mathfrak{CM}}
\newcommand{\mOp}{\mathcal{MO}p}
\newcommand{\MD}{\mathfrak{M}^{\D}}
\newcommand{\F}{\mathcal{F}}
\newcommand{\Ff}{\mathbb{F}}
\newcommand{\nthord}{^{(n)}}
\newcommand{\Aut}{\underline{\operatorname{Aut}}}
\newcommand{\Gr}{{\mathfrak{Gr}}}
\newcommand{\GR}{\on{Gr}}
\newcommand{\GRo}{{\mathfrak{GR}^{\circ}}}
\newcommand{\Fr}{\operatorname{Fr}}
\newcommand{\GL}{\operatorname{GL}}
\newcommand{\gl}{\mathfrak{gl}}
\newcommand{\SL}{\operatorname{SL}}
\newcommand{\ff}{\footnote}
\newcommand{\ot}{\otimes}
\def\Ext{\operatorname {Ext}}
\def\Hom{\operatorname {Hom}}
\def\Ind{\operatorname {Ind}}
\def\bbZ{{\mathbb Z}}

\newcommand{\AutO}{\on{Aut}\Oo}
\newcommand{\Der}{{\on{Der}\,}}
\newcommand{\DerO}{{\on{Der}\Oo}}
\newcommand{\AutK}{\on{Aut}\K}

\newcommand{\nc}{\newcommand}
\nc{\ol}{\overline}
\nc{\cont}{\on{cont}}
\nc{\rmod}{\on{mod}}
\nc{\Mtil}{\widetilde{M}}
\nc{\wb}{\overline}
\nc{\wt}{\widetilde}
\nc{\wh}{\widehat}
\nc{\sm}{\setminus}
\nc{\mc}{\mathcal}
\nc{\mbb}{\mathbb}
\nc{\Mbar}{\wb{M}}
\nc{\Nbar}{\wb{N}}
\nc{\Mhat}{\wh{M}}
\nc{\pihat}{\wh{\pi}}
\nc{\JYX}{\cJ_{Y\leftarrow X}}
\nc{\phitil}{\wt{\phi}}
\nc{\Qbar}{\wb{Q}}
\nc{\DYX}{\D_{Y\leftarrow X}}
\nc{\DXY}{\D_{X\to Y}}
\nc{\dR}{\stackrel{\bbL}{\underset{\D_X}{\ot}}}
\nc{\Winfi}{\cW_{1+\infty}}
\nc{\K}{{\mc K}}
\nc{\Kx}{{\mc K}^{\times}}
\nc{\Ox}{{\mc O}^{\times}}
\nc{\unit}{{\bf \on{unit}}}
\nc{\boxt}{\boxtimes}
\nc{\xarr}{\stackrel{\rightarrow}{x}}

\nc{\Gamx}{\Gamma_-^{\times}}
\nc{\Gap}{\Gamma_+}
\nc{\Emx}{\cE_-^{\times}}

\nc{\Enat}{E^{\natural}}
\nc{\Enatbar}{\ol{E}^{\natural}}
\nc{\Of}{{\mathbb O}}
\nc{\obar}{\wb{o}}
\newcommand{\SShv}{\on{SShv}}

\nc{\Ga}{\G_a}

\nc{\Gm}{\G_m}
\nc{\Gabar}{\wb{\G}_a}
\nc{\Gmbar}{\wb{\G}_m}

\title{From Solitons to Many--Body Systems}
\author{David Ben-Zvi}
\address{Department of Mathematics\\University of Texas\\Austin, TX 78712-0257}
\email{benzvi@math.utexas.edu}
\author{Thomas Nevins}
\address{Department of Mathematics\\University of Michigan\\Ann Arbor, MI 48109-1109}
\email{nevins@umich.edu}

\begin{abstract}
We present a bridge between the KP soliton equations and the
Calogero--Moser many-body systems through noncommutative algebraic
geometry. The Calogero-Moser systems have a natural geometric
interpretation as flows on spaces of spectral curves on a ruled
surface. We explain how the meromorphic solutions of the KP hierarchy
have an interpretation via a {\em noncommutative} ruled
surface. Namely, we identify KP Lax operators with vector bundles on
quantized cotangent spaces (formulated technically in terms of
$\D$-modules). A geometric duality (a variant of the Fourier--Mukai
transform) then identifies the parameter space for such vector bundles
with that for the spectral curves and sends the KP flows to the
Calogero--Moser flows.  It follows that the motion and collisions of
the poles of the rational, trigonometric and elliptic solutions of the
KP hierarchy, as well as of its multicomponent analogs, are governed
by the (spin) Calogero--Moser systems on cuspidal, nodal and smooth
genus one curves. This provides geometric explanations and
generalizations of results of Airault--McKean--Moser, Krichever and
Wilson. The present paper is an overview of work to
appear in \cite{solitons}.
% {\it To appear, Proceedings of MSRI
% Special Semester on Algebraic Stacks.}
\end{abstract}

\maketitle

\section{Introduction} 
Our purpose in this paper is to introduce a geometric viewpoint
(developed in detail in \cite{solitons}) on a much-explored, puzzling
phenomenon of the theory of integrable systems: the description of the
motion of poles of meromorphic solutions of soliton equations by
simple many-body systems. Our viewpoint provides a conceptual
framework for this phenomenon within noncommutative algebraic
geometry. As an application we obtain direct geometric proofs of
results that complete the work of Airault--McKean--Moser \cite{AMM},
Krichever \cite{Kr1,Kr2} and Wilson \cite{Wilson CM} on the motions
and collisions of poles of meromorphic KP solutions.

\subsection{Historical Overview}
Before reviewing the Calogero--Moser system (CM), the
Kadomtsev--Petviashvili hierarchy (KP), and their relation, we sketch
an incomplete historical overview of the problem---see the review
articles \cite{Benn,BAMS} for more complete history and
bibliography. In the seminal work \cite{AMM}, Airault, McKean and
Moser wrote down rational, trigonometric and elliptic solutions of the
Korteweg--deVries equation and discovered that the motion of their
poles is governed by the Calogero--Moser classical many-body systems
of particles on the line, cylinder and torus (respectively) with
inverse square potentials. Krichever \cite{Kr1,Kr2} and the
Chudnovskys \cite{CC} extended this correspondence to the meromorphic
solutions of the KP equation, where it becomes an isomorphism between
the phase spaces of generic rational (decaying at infinity),
trigonometric and elliptic KP solitons and the corresponding
Calogero--Moser systems. 

Krichever derived this result (in the
elliptic case) from a relation between an auxiliary linear problem associated to a KP potential and the Calogero--Moser systems.  
More precisely, to an elliptic 
KP potential one associates a non-stationary Schr\"odinger operator
with elliptic potential, and Krichever proved that this auxiliary 
Schr\"odinger operator has meromorphic solutions if and only if
the poles of the
potential move as particles in the elliptic CM system. 
In particular, this
showed that the elliptic CM system can be written in terms of spectral
curves (using a Lax operator with spectral parameter), and showed that
the generic elliptic solutions of KP are {\em finite gap solutions}---they 
come from applying Krichever's general geometric construction
of solutions to KP to these spectral curves, known as {\em tangential
covers} of the elliptic curve. A detailed algebro-geometric study of
tangential covers was undertaken by Treibich and Verdier \cite{TV,TV2},
leading to a complete classification of elliptic solutions of the KdV
equation. This geometric description of the elliptic CM systems may be
used to identify them with a meromorphic version of the Hitchin system
\cite{GorNe,Nekrasov,DW,DM}---for the rational and trigonometric
systems, the corresponding Hitchin systems live on cuspidal or nodal
(rather than smooth) genus one curves \cite{Nekrasov} (see also
\cite{Kr5} for a different point of view on this identification,
closer in spirit to the current work, as part of the theory of
meromorphic Lax operators on curves).

Other work on the KP/CM correspondence includes generalizations to
some cases of multicomponent KP hierarchy and spin Calogero--Moser
systems \cite{BBKT,T matrix}, an analog for difference equations
relating the 2D Toda hierarchy and Ruijsenaars-Schneider systems
\cite{KrZab,T difference,Kr3,Kr4}, and extensions to other related
systems (see e.g. \cite{Kr3,Kr4,BAMS,BB,Braden}). The KP/CM
correspondence is also applied in the study of the bispectral
phenomenon \cite{DG,Wilson bispectral,Wsurv1, Wilson CM,Wsurv2,BW
automorphisms,Wsurv3} and relates to Seiberg--Witten integrability of
supersymmetric Yang--Mills theory (see \cite{Whitham} for a collection
of reviews).

\subsection{Some Recent Developments.}
The KP/CM correspondence in the rational case was greatly deepened by
Wilson \cite{Wilson CM} (see \cite{Wsurv3} for a review), who extended
it to allow collisions of particles.  More precisely, the rational
Calogero--Moser phase space possesses a natural completion that is
constructed as a space of pairs of matrices whose commutator lies in a
particular conjugacy class. In this description, the positions of the
particles correspond to the eigenvalues of one of the matrices, which
are now allowed to coincide.  Wilson had (see \cite{Wilson
bispectral}) identified the completed phase space of the rational KP
hierarchy with an {\em adelic Grassmannian} that parametrizes certain
subspaces of $\C[x]$. In \cite{Wilson CM}, Wilson gives an explicit
formula that defines a point of the adelic Grassmannian from a CM pair
of matrices and then proves by direct calculations that this map
extends continuously to the completed phase spaces and takes the CM
flows to the KP flows.  (In the generic rational case, Shiota
\cite{Sh} has previously extended the KP/CM correspondence to all the
higher flows of the KP hierarchy and the higher CM hamiltonians,
establishing a bijection between generic rational solutions of the KP
{\em hierarchy} and the Calogero--Moser hierarchy; see also the
related work \cite{AKV}).

Wilson's adelic
Grassmannian appears independently in the work of Cannings and Holland
\cite{CH ideals} classifying (right) ideals in the Weyl algebra of
differential operators on the affine line, indicating that there might
be an interesting relationship between the KP/CM correspondence and
objects of noncommutative algebra/geometry. In \cite{BW automorphisms}
(see \cite{Wsurv2} for a review), Berest and Wilson show that the
decomposition of the rational KP phase space into the union of
$n$-particle CM spaces can be described as the orbit decomposition
under the action of the group of automorphisms of the Weyl algebra,
``half'' of which is matched up with the KP flows. This decription is
then used to identify Wilson's bispectral involution \cite{Wilson
bispectral} (see \cite{Wsurv1} for a review) on the rational KP
solutions with the Fourier transform on the Weyl algebra. 

In the
papers \cite{BW ideals} by Berest and Wilson and \cite{BGK1} by
Baranovsky, Ginzburg and Kuznetsov, a direct relation between the
Calogero--Moser pairs of matrices and the classification of right ideals in
the Weyl algebra is presented through a calculation in noncommutative
geometry. Namely, by interpreting the ideals as sheaves on a
noncommutative projective plane, one can apply the techniques of
Koszul duality and the Beilinson spectral sequence to classify
ideals (up to isomorphism as modules) by cohomological data (monads), which turns out
to reproduce precisely the Calogero--Moser matrices. 

However, a direct
relation between KP and ideals in the Weyl algebra (or noncommutative
geometry) was missing, as were a conceptual explanation of the relation
between KP and Calogero--Moser systems and an extension of the
correspondence to completed phase spaces in the trigonometric,
elliptic and multicomponent cases. These goals are achieved in
\cite{solitons}, and described in the present paper.

\subsection{The Current Work.}
As we have already indicated, the present work is devoted to to
introducing a geometric approach to the KP/CM correspondence that is
developed in detail in \cite{solitons}.  Accordingly, we begin with a
review (Section \ref{CM section}) of the (rational, trigonometric and
elliptic) Calogero--Moser systems and their concrete construction (by
pairs of matrices, in the rational case).  We then describe the
formulation \cite{TV2}
 of the Calogero--Moser systems as flows on the space of
pairs $(\Sigma,\cL)$ of a line bundle $\cL$ on a curve $\Sigma$ (the
{\em spectral curve}), embedded in a ruled surface $\Enat$ over a
(cuspidal, nodal or smooth) genus 1 curve $E$.
 To pass between this
geometric picture and the concrete Calogero--Moser particles, we
recall the identification of both with a special case of the
meromorphic Hitchin system on $E$.

In Section \ref{KP} we discuss the KP hierarchy, its Lax formulation
in terms of formal microdifferential operators, and Sato's
reformulation of KP Lax operators in terms of $\cD$-modules. Sato
identifies Lax operators with an open subset (the big cell
$\GR^{\circ}$) of an infinite-dimensional Grassmannian $\GR$.  We
introduce a coordinate-free reformulation of {\em meromorphic} KP Lax
operators, the {\em micro-opers}, by generalizing the $\D$-modules
appearing in Sato's theory. We then extend (Theorem \ref{microopers and
GR}) Sato's identification between the big cell of the Grassmannian
and the set of Lax operators on the disc to an identification of the entire
Grassmannian with the space of micro-opers on the disc. Micro-opers
on a general curve $X$ define meromorphic Lax operators on $X$, whose
poles correspond to the points where the ``local data'' lie outside
the big cell of the Sato Grassmannian. Thus micro-opers are
well-suited to describing the poles of Lax operators and their
collisions. The formulation of the KP hierarchy in terms of (regular)
micro-opers is closely analogous to the Drinfeld--Sokolov
formulation of the KdV hierarchies in terms of connections---indeed,
micro-opers are the KP analogues of the {\em opers} of
Beilinson--Drinfeld (\cite{Hecke}) (or more precisely of the {\em
affine opers} of \cite{BF}).

It is at this stage that noncommutative algebraic geometry enters the
picture: we explain in Section \ref{D-bundles} how a micro-oper on a
curve $X$ may be interpreted as a line bundle or rank one torsion free
sheaf on a noncommutative ruled surface, (the completion of) the {\em
quantized cotangent bundle} $T^*_{\hbar}X$ of $X$, equipped with some data along a
divisor ``at infinity''.  In this setting micro-differential operators
arise geometrically as Laurent expansions of functions on the
noncommutative surface along this divisor. Moreover, there is a simple
geometric description of the KP flows directly on the space of
micro-opers: in terms of the noncommutative ruled surface, the flows
act as modifications of bundles along the divisor by ``changing the
transition function.''

In Section \ref{Fourier} we explain how a geometric Fourier transform
provides the KP/CM correspondence.  Laumon \cite{La} and Rothstein
\cite{Ro2} enhanced the original Fourier--Mukai transform for abelian
varieties to a geometric integral transform that takes $\D$-modules on
an abelian variety $A$ to quasicoherent sheaves on a bundle
$\widehat{A}^{\natural}$ over the dual abelian variety $\widehat{A}$.
This enhanced Fourier--Mukai transform was then used in work of
Nakayashiki \cite{N1,N2} and Rothstein \cite{Ro1,Ro2} to describe the
algebro-geometric solutions of KP. We extend the enhanced
Fourier--Mukai transform to the case of singular genus one curves (see
Theorem \ref{thm5.2}) and use it to construct an isomorphism between
moduli spaces of vector bundles on noncommutative ruled surfaces over
genus one curves and moduli spaces of spectral sheaves on a
(commutative) ruled surface.  This gives a new proof of the
classification \cite{BW ideals, BGK1} of ideals in the Weyl algebra by
quiver data, the Calogero--Moser pairs of matrices (this follows from
the rank one, rational case of our result).  The same technique in a
degenerate case gives a new proof of the ADHM classification of framed
torsion-free sheaves on $\Pp^2$ by quiver data.

Our main result, described in Section \ref{final}, states that the
Fourier transform precisely identifies micro-opers on a genus one
curve with the spectral data of the CM system, and the KP flows on the
former with the CM system on the latter:

\begin{thm}[\cite{solitons}]\label{microops and CM} Let $E$ denote an irreducible genus one curve.
The extended Fourier--Mukai transform induces an isomorphism of the
moduli spaces of micro-opers on $X$ and of Calogero-Moser spectral
data that identifies the KP flows with the Calogero--Moser flows. The poles of the
micro-opers correspond to the positions of Calogero--Moser particles,
extending the bijection of \cite{Kr1,Kr2} in the generic case and that
of \cite{Wilson CM} in the completed rational case.
\end{thm} 

The KP solutions with {\em generic} singularities, which are
identified with configurations of CM particles with
{\em distinct} positions in \cite{Kr1,Kr2},
correspond to micro-opers which only hit the codimension one strata
of the Sato Grassmannian. The collisions of CM particles correspond
to passing to more interesting singularities of the micro-opers, i.e.
deeper strata of the Grassmannian. In particular, we obtain that
the correspondence extends to completed phase spaces for the rational,
trigonometric and elliptic systems, generalizing Wilson's result
\cite{Wilson CM} in the rational one-component case.
(Here Wilson's ad\`elic Grassmannian appears naturally in its
Cannings--Holland interpretation, as a parameter space of
$\D$--modules.)

The same technique, applied to higher rank vector bundles or
micro-opers, identifies the completed phase spaces of the
multicomponent generalizations of KP with those for the spin
generalizations of the Calogero--Moser systems, as is described in
\cite{solitons}.

We conclude in Section \ref{further aspects} with brief comments on
further work on the relation with $\cW_{1+\infty}$ vertex algebras
\cite{W}, the bispectral involution of Wilson, and generalization of
the results of \cite{solitons} to the difference analog
(Toda/Ruijsenaars correspondence) \cite{solitons2}.

\subsection{Calogero--Moser Particles as 
Points on a Noncommutative Surface}  We have already discussed the role of noncommutative geometry in a geometric formulation of meromorphic solutions of KP.  However, Theorem \ref{microops and CM} also
has an interesting noncommutative geometric interpretation that is close in
spirit to familiar commutative constructions.  

Namely,
recall that the Hilbert scheme of points on the cotangent bundle $T^*X$ of a
curve $X$ is isomorphic to (a component of) the moduli space of rank one
torsion-free sheaves on the smooth projective completion 
$\overline{T^*X} = T^*X \cup X$ 
that are {\em framed} (trivialized) along the divisor at
infinity.   When one deforms $T^*X$ to the quantized cotangent bundle
 $T^*_{\hbar}X$ of 
$X$ it is this latter description of the Hilbert scheme that deforms well;
consequently one takes, as the Hilbert scheme of points on $T^*_{\hbar}X$, 
the moduli space of framed rank one torsion-free sheaves on the 
projective completion $\overline{T^*_{\hbar}X}$.  
Note that every micro-oper on $X$ gives such a
 framed rank one torsion-free
sheaf on $\overline{T^*_{\hbar}X}$.
On the other hand, the natural phase space for a many-body system
(such as the CM system) on an elliptic curve $E$ is the configuration
space for collections of
 points on the cotangent bundle $T^*E$ with distinct
$E$-coordinates. We may seek to complete the phase space (and
Hamiltonians) so as to allow collisions of particles. An obvious
candidate for a completed phase space would be the Hilbert scheme of
points on the cotangent bundle. The natural dynamical system on the
Hilbert scheme of points on $T^*E$, however, is the trivial many-body
system, with zero potential between distinct particles.

In this language our result says that the Calogero--Moser system
extends to the Hilbert scheme of points on the {\em quantized}
cotangent bundle of the curve. While the quantized cotangent bundle
does not have any points in the ordinary sense (since $\D_E$ has no
finite-dimensional modules), the space of micro-opers provides a
natural candidate for this Hilbert scheme. We thus show that the
natural dynamics of ``points on the quantized cotangent bundle'' of
$E$ (that is, of micro-opers on $E$) is the KP system, and identify this
precisely with the Calogero-Moser system. Moreover the positions of
the Calogero--Moser particles are easily read off from the micro-oper
(as the $E$-coordinates of the putative points in $T^*_{\hbar}E$), 
and are identified with the poles of KP Lax operators.  

From
this point of view, both KP and Calogero--Moser systems naturally live
on the quantized cotangent bundle---the interpretation using spectral
curves via the Fourier--Mukai transform is then a tool to
 describe the relevant moduli spaces and check that we
are indeed getting the correct Hamiltonian system.

\subsection{Acknowledgements.}
We would like to thank Ron Donagi, Dennis Gaitsgory, Victor Ginzburg
and Tony Pantev for very helpful conversations. We are especially
grateful to Edward Frenkel and Mitch Rothstein for generously sharing
their understanding of the geometry of soliton equations. Both authors
were supported during the course of this work by NSF postdoctoral
fellowships, and by MSRI postdoctoral fellowships as part of the
Spring 2002 programs ``Algebraic Stacks, Intersection Theory and
Non-Abelian Hodge Theory'' and ``Infinite-Dimensional Algebras and
Mathematical Physics.''

\section{Calogero--Moser Systems}\label{CM section}
We begin by reviewing the complexified Calogero-Moser systems
associated with the general linear group $GL_n$, following
\cite{KKS,Wilson CM,Nekrasov}.

\subsection{The CM System}
Recall that connected one-dimensional complex algebraic groups $\G$
fall into three classes: the additive group $\Ga=\C$, the
multiplicative group $\Gm=\Cx$, and the one-parameter family of
elliptic curves $E$. We denote the identity element of each of these
groups by $o$. The complexified Calogero--Moser systems are completely
integrable Hamiltonian systems describing a collection of $n$
identical particles on one of these groups $\G$. Thus the phase space
of the Calogero--Moser system on the curve $\G$ is
$T^*\big(\G^{(n)}\sm\bigcup \{\text{all diagonals}\}\big)$, 
the configurations of $n$ distinct unlabelled points
$q_i\in \G$ with momenta $p_i\in\C$. They are described, in terms of
coordinates on the complex line (the universal cover of $\G$) by the
Hamiltonians
$$H=\frac{1}{2}\sum_{i=1}^n p_i^2 + \sum_{i<j} U(q_i-q_j),$$ with
potential functions with a single second order pole at the origin of
$\G$, that is
\begin{equation*}
{\bf Rational: }U(q)=\frac{1}{q^2}\hskip.2in {\bf Trigonometric:
}U(q)=\frac{1}{\sin^2(q)}
\hskip.2in {\bf Elliptic: }U(q)=\wp(q)
\end{equation*}
where $\wp(q)$ is the Weierstrass $\wp$-function attached to the
elliptic curve $E$.

\subsubsection{}\label{Kostant CM}
The rational and trigonometric Calogero-Moser systems can be described
concretely by Hamiltonian reduction for the group $GL_n$ (\cite{KKS}).
Let $\Of\subset\gl_n$ denote the conjugacy class of traceless matrices
of the form $\on{Id}+A$ with $A$ a rank one matrix (which we consider
as a coadjoint orbit in $\gl_n^*$). It is convenient to realize $\Of$
as the orbit of either of the matrices
\begin{equation}\label{CM orbit}\left(\begin{array}{ccccc}
0&1&1&\cdots&1\\ 1&0&1&\cdots&1\\ 1&1&0&\cdots&1\\
\vdots&\vdots&\hdots&\ddots&\vdots\\ 1&1&1&\cdots&0
\end{array}\right) \;,\;
\left(\begin{array}{ccccc}
1-n&0&0&\cdots&0\\
0&1&0&\cdots&0\\
0&0&1&\cdots&0\\
\vdots&\vdots&\hdots&\ddots&\vdots\\
0&0&0&\cdots&1
\end{array}\right).
\end{equation}
Then the Hamiltonian reduction of $T^*\gl_n$ by $GL_n$
with the moment condition $\Of$ is the space
$$\CM_n= \lbrace (X,Y)\in T^*\gl_n \;\big|\; [X,Y] \in
\Of\rbrace/GL_n.$$ As proven in \cite{Wilson CM}, this space is a
smooth irreducible affine variety of dimension $2n$. This variety
comes equipped with the Hamiltonian function
$H=\frac{1}{2}\on{tr}Y^2$, and all other invariant polynomials of $Y$
give Hamiltonians in involution with $H$, whence one deduces the
complete integrability of the Calogero-Moser system by a dimension
count. On the open subset where $X$ has distinct eigenvalues, we may
write coordinates $(q_i,p_i)$ on (a finite cover of) $\CM_n$, using
the first matrix representative of $\Of$ in \ref{CM orbit}:
$$X=\left(\begin{array}{ccccc} q_1&0&0&\cdots&0\\ 0&q_2&0&\cdots&0\\
0&0&q_3&\cdots&0\\
\vdots&\vdots&\hdots&\ddots&\vdots\\
0&0&0&\cdots&q_n
\end{array}\right), \hskip.2in
Y=\left(\begin{array}{ccccc}
p_1&\frac{1}{q_1-q_2}&\frac{1}{q_1-q_3}&\cdots&\frac{1}{q_1-q_n}\\
\frac{1}{q_2-q_1}&p_2&\frac{1}{q_2-q_3}&\cdots&\frac{1}{q_2-q_n}\\
\frac{1}{q_3-q_1}&\frac{1}{q_3-q_2}&p_3&\cdots&\frac{1}{q_3-q_n}\\
\vdots&\vdots&\hdots&\ddots&\vdots\\
\frac{1}{q_n-q_1}&\frac{1}{q_n-q_2}&\frac{1}{q_n-q_3}&\cdots&p_n
\end{array}\right).
$$
It is easy to see that the Hamiltonian $H$ in these coordinates
recovers the rational Calogero--Moser Hamiltonian above. Thus $\CM_n$
provides a completion of the phase space of the rational
Calogero--Moser system, in which we allow the points $q_i$ (the
eigenvalues of $X$) to collide.  The trigonometric Calogero-Moser
system can be described similarly from reduction of $T^*GL_n$:
we consider conjugacy classes of pairs of matrices $(X,Y)$ with $X$
invertible and $XYX\inv-Y\in\Of$. As before, we obtain explicit
coordinates on the locus where $X$ is diagonalizable and the commuting
Hamiltonians from the invariant polynomials of $Y$. In the next
section we will see how the Calogero--Moser rank one matrices arise
naturally from the consideration of spectral curves over a singular
elliptic curve.

\subsection{Calogero--Moser Spectral Curves.}\label{cubic curves}\label{Calogero and Hitchin}
In this section we present the geometric construction of the
Calogero--Moser systems in terms of spectral curves on an algebraic
surface (see \cite{TV,TV2,GorNe,DW,DM} for the elliptic case, and
\cite{Nekrasov} for the rational and trigonometric cases). The
connection with the explicit description of the previous section is
explained in Section \ref{CM Hitchin} by interpreting both as Hitchin
systems.

A general paradigm for constructing integrable systems (see \cite{DM} for
a discussion and references) involves fixing a
symplectic algebraic surface and a family of curves (specifically, a
linear series) on this surface. The phase space of the associated
integrable system is the space of pairs $(\Sigma,\cL)$ of a curve
$\Sigma$ in the family and a line bundle (or, for singular $\Sigma$, a
rank one torsion--free sheaf) $\cL$ on $\Sigma$ (i.e. the family of
generalized Jacobians of the curves in the linear series). The
Hamiltonians of the system are given by coordinates on the space of
curves $\Sigma$, and the flows of the system are linear flows on the
Jacobians, modifying $\cL$ while fixing $\Sigma$.

In the case of the rational, trigonometric and elliptic
Calogero--Moser systems, the relevant symplectic surface $\Enat$ is a
special affine bundle over a cuspidal, nodal or smooth elliptic curve
$E$, respectively.  Recall that an irreducible reduced algebraic curve
$E$ of arithmetic genus one, i.e. a Weierstrass cubic curve, falls
into one of three classes:
\begin{enumerate}
\item[$\bullet$] {\bf Elliptic:} $E$ is a smooth elliptic curve (in
particular a group), and may be described by an equation of the form
$y^2=x^3+a x+b$.
\item[$\bullet$]{\bf Trigonometric:} $E=\Gmbar$ is a nodal cubic, and
is isomorphic to the curve $y^2=x^2(x-1)$.  Its normalization
$\pline\to E$ identifies two points $0$ and $\infty$ to a node on $E$,
and defines a group structure $\Gm\cong E^{sm}\subset E$ on the smooth
locus.
\item[$\bullet$]{\bf Rational:} $E=\Gabar$ is a cuspidal cubic, and is
isomorphic to the curve $y^2=x^3$. Its normalization $\pline\to E$
collapses $2\cdot\infty$ to a cusp on $E$, and defines a group
structure $\Ga\cong E^{sm}\subset E$ on the smooth locus.
\end{enumerate}
In each of the cases we have canonical isomorphisms $E^{sm}=\G\cong
\on{Pic}^0 E$ of the corresponding group $\G$ with the Jacobian of
$E$, and of $E$ itself with the compactified Jacobian of $E$, the
moduli of torsion free sheaves of rank one and degree zero on $E$
(trivialized at the identity).  
\footnote{In particular there is a universal
sheaf $\cP\to E\times E$, the {\em Poincar\'e sheaf}.}

An elliptic curve $E$ carries a unique nontrivial rank 1 affine bundle
$\Enat\to E$. Let $\on{At}$ denote the {\em Atiyah bundle} on $E$, the
unique nontrivial extension of $\Oo_E$ by itself. The (Stein)
algebraic surface $\Enat$ may be identified with the complement of the
unique section $E_\infty=\Pp(\Oo)\cong E$ of the projectivization of
the Atiyah bundle, $\Enat\subset\Enatbar=\Pp(\on{At})\subset
E_\infty$. $\Enat$ is also naturally identified with the moduli space
of line bundles with a holomorphic connection on $E$.
\footnote{In particular,
the pullback of the Poincar\'e bundle to $E\times \Enat$ carries a
canonical connection relative to the second factor.}  Finally, $\Enat$
may also be identified as a twisted cotangent bundle of $E$, namely as
the affine bundle of connections on the line bundle $\Oo(o)$ on
$E$. In particular it inherits a symplectic form.  

For a singular cubic curve, one defines an affine bundle $\Enat\to E$
with analogous properties, taking account of the behavior at the
singularity (effectively replacing the tangent sheaf of the curve by
its trivial line subbundle generated by the invariant vector field on
the group $\G$)---see
\cite{solitons} for a detailed account. In each case, the projective
bundle $\Enatbar\to E$ has a unique section $E_\infty$, whose
complement is the affine bundle $\Enat$.

The open (Stein) surface $\Enat$ does not contain any complete curves,
so that any curve in $\Enatbar$ must intersect the divisor
$E_\infty$. So it is natural to consider the simplest linear series on
$\Enat$, consisting of curves with only one point of intersection with
$E_\infty$ (taken for convenience to be the origin $o\in E_\infty$),
and degree $n$ over $E$ for some positive integer $n$.  This linear
series (denoted $|n\cdot E_\infty+F|$) on $\Enatbar$ is studied in
detail in \cite{TV}.  In particular, it is equivalent for a curve
$\Sigma\in |n\cdot E_\infty+F|$ to be irreducible and to have a unique
point of intersection with the section $E$ at infinity (which is then
automatically transversal).  Moreover the collection of such curves
$\Sigma$ is the complement of a hyperplane in the projective space
$|n\cdot E_\infty+F|$. There is also a bijection between such
$\Sigma\hookrightarrow \Enatbar\to E$ and {\em tangential covers} of
$E$, corresponding to maps $E\to\on{Jac}\Sigma$ which are tangent at
$o$ to the Abel--Jacobi line of $\Sigma$ at the distinguished point
$\Sigma\cap E_\infty$.

Thus we consider the integrable system, the {\em completed
Calogero--Moser system}, with phase space given by line bundles, and
more generally rank one torsion--free sheaves, supported on such
curves:

\begin{defn} \label{CM sheaf}
Let $\Oo_o$ denote the skyscraper sheaf at $o$.  A {\em CM spectral
sheaf} is a coherent sheaf $\cL$ on $\Enatbar$ of pure dimension one,
equipped with an isomorphism $\phi:\cL|_{E_{\infty}}\to \Oo_o$, and with
first Chern class $n\cdot E_\infty+F$ (where $F$ denotes the fiber
over $o$).\footnote{If $E$ is singular we add the technical condition
that $M$ is locally free over the singular point.} The space of CM
spectral sheaves is denoted by $\CM_n(E)$.
\end{defn}

The Hamiltonian flows of the completed Calogero--Moser system preserve
 the underlying curve $\Sigma$ while modifying the spectral sheaf
 (linearly) along the generalized Jacobian of $\Sigma$. We may
 generate all such flows by modifying $\cL$ only at the distinguished
 point $o$ at infinity (as in \cite{BL,BF}). Let $(\cL,\Sigma)\in
 \CM_n(E)$ be a Calogero-Moser spectral sheaf, and let
 $\K_{\Sigma,\infty}$ denote the field of Laurent series along
 $\Sigma$ at $\Sigma\cap E_\infty$. Then there is a canonical
 identification
 $\on{End}(\cL\ot\K_{\Sigma,\infty})=\K_{\Sigma,\infty}$. It easily
 follows that the abelian Lie algebra $\K_{\Sigma,\infty}$ acts
 (linearly and formally transitively) on the compactified Picard of
 $\Sigma$ by formally changing the transition function of $\cL$ at
 $\Sigma\cap E_\infty$.

\subsection{CM and Hitchin}\label{CM Hitchin}
In this section we sketch the description of the Calogero--Moser systems
 as Hitchin systems following \cite{Nekrasov,DM}.

Let $Bun_n(E,o)$ denote the moduli space of rank $n$ semistable vector
bundles of degree zero on a cubic curve $E$, equipped with a
trivialization of the fiber at the identity; for singular $E$, we
impose the (open) condition that the bundle have trivial pullback to
the normalization $\pline$. The cotangent fiber $T^*Bun_n(E,o)|_V$ at
a bundle $V$ consists of pairs $(V,\eta)$, where $\eta$ is a 
meromorphic Higgs field
$\eta\in\Gamma(\on{End}V(o))$ on $V$ with only a simple pole at
$o$. We may perform Hamiltonian reduction of $T^*Bun_n(E,o)$ with
respect to the action of $\GL_n$ changing the trivialization at $o$,
with moment map taking values in an arbitrary coadjoint orbit in
$\gl_n$.  Thus we define the Calogero--Moser--Hitchin space 
by
\begin{equation}\label{Hitchin CM}
T^*Bun_n(E,o)\,/\!\!/_{\Of}\, \GL_n  = \big\{ (V,\eta)\; :\;\on{res}_o(\eta)\in\Of\big\}/\GL_n,
\end{equation}
that is, as the space of Higgs bundles on $E$ having residue at $o$ in
the Calogero--Moser orbit $\Of$. The description as a reduction endows
the Calogero--Moser--Hitchin space with a symplectic form and with $n$
algebraically independent Poisson commuting Hamiltonians given by the
invariant polynomials in $\eta$ (the Hitchin hamiltonians).

In the rational and trigonometric cases one has concrete descriptions
of the moduli spaces: since our bundles have trivial pullback to the
normalization $\pline$, they are completely described by the descent
data from $\pline$ to $E$. This descent data in the nodal case is the
identification of the two fibers over the inverse image of the node,
hence $Bun_n(\Gmbar,o)=\GL_n$. In the cuspidal case these two points
are infinitesimally nearby, and the descent data becomes a
``connection matrix'' identifying these two nearby fibers---thus we
have $Bun_n(\Gabar,o)=\gl_n$.  It follows immediately that the
Calogero--Moser--Hitchin space for $\Gabar$ is isomorphic with the
completed rational Calogero--Moser phase space $\CM_n$ introduced in
Section \ref{Kostant CM}, and similarly in the trigonometric case.  In
the elliptic case we may recover the elliptic Calogero--Moser particles and
Hamiltonians as
follows. There is a dense open subvariety of $Bun_n(E,o)$
that is identified with the configuration space of $n$
distinct points $q_i$ on $E$, via the assignment
\begin{equation}\label{elliptic particles}
\{q_i\}\mapsto \oplus
\Oo(q_i-o).
\end{equation} 
Writing the Hitchin Hamiltonian in these coordinates
gives the elliptic Calogero-Moser flow on the $q_i$.

Hitchin's integrable system on a curve $X$ is naturally described in
terms of spectral curves, this time embedded in the symplectic surface
$T^*X$. In our case, the minimal (rank one) condition for the matrix
$A$ in the coadjoint orbit condition on the Calogero--Moser matrices
$[X,Y]=\on{Id}+A$ matches up precisely with the minimal (one-point)
intersection condition on the Calogero--Moser spectral curves, yielding the
following:

\begin{prop}\label{CM and SShv} The $n$ particle Calogero--Moser phase space 
on $E$ is canonically identified with an open subset of the space
 $\CM_n(E)$ of 
Calogero--Moser spectral sheaves with support of degree $n$ over $E$.
Under this isomorphism, the CM flows are identified with linear flows
along the generalized Jacobians of CM spectral curves.
\end{prop}

\begin{remark} See \cite{DW} for a geometric description of this 
identification via the CM Hitchin system.  One begins by
describing the meromorphic Hitchin system in terms of spectral curves
$\wt{\Sigma}$ in the total space of the line bundle
$\Omega_E(0)=\Oo_E(0)$ of differentials with simple pole at $0$. One
then obtains a CM spectral curve through a birational transformation,
blowing up the point in the fiber over $0$ corresponding to the
eigenvalue $1$, and blowing down the proper transform of the fiber.
\end{remark}

\section{The KP Hierarchy}\label{KP}

In this section we briefly review the definition of the KP hierarchy
and its interpretation via flows on the Sato Grassmannian.  General
references for this section are \cite{Mu}, \cite{Sa} and \cite{SW}.

\subsection{Introducing KP}
Let $\cE$ denote the algebra of formal microdifferential operators (or
pseudodifferential symbols) with coefficients in $\C[[t]]$. An element
of $\cE$ is a Laurent series
\begin{equation}\label{microdiff}
M=\sum_{N\ll \infty}a_N\del^N\hskip.3in a_i\in\C[[t]]
\end{equation} in the
formal inverse $\del\inv$ of the derivation $\del=\del_t$ of
$\C[[t]]$. The composition in $\cE$ is determined by the Leibniz rule,
\begin{equation}\label{Leibniz}
\del^n \cdot f=\sum_{i\geq 0}\binom{n}{i} f^{(i)}\del^{n-i},
\end{equation}
where $\binom{n}{i}$ is defined for $n<0$ by taking
\bd
\binom{n}{i} = \frac{n(n-1)\cdots (n-i+1)}{i(i-1)\cdots 2\cdot 1}.
\ed 
The
ring $\D=\C[[t]]\langle\del \rangle$ of differential operators is a
subring of $\cE$. We also have a commutative subring $\Gamma=\C((\del\inv))$ of
{\em constant coefficient} microdifferential operators.

Consider a
microdifferential operator of the form
\bd
L=\partial+u_1\partial\inv+u_2\partial^{-2}+\cdots\in\cE,
\ed
which is called a KP {\em Lax operator}.  
The space $\cL$ of such operators is an infinite-dimensional
affine space (with coordinates the coefficients of the $u_i$).  The KP hierarchy
is the collection of compatible evolution equations on a Lax operator $L$
defined as follows:
\begin{equation}\label{KP flows}
\frac{\partial L}{\partial t_n}=[L,(L^n)_+],
\end{equation} where
$(M)_+=\sum_{N\geq 0} a_N \del^N\in\D\subset\cE$ denotes the
differential part of a microdifferential operator $M$ as in 
\eqref{microdiff}.  That is, we let the operator
$L=L(t,t_1,t_2,\dots)$ depend on the infinitely many time variables
$t_n$ and then require that the dependence of $L$ on $t_n$ (i.e. its
``evolution along the $n$th time'') satisfies \eqref{KP flows}. 

 This may be
immediately reinterpreted in terms 
of a collection of vector fields on $\cL$: we define the $n$th 
vector field on the affine space $\cL$ by taking its value at $L$ 
to be the commutator $[L,(L_n)_+]$.  A solution $L$ of the equations
\eqref{KP flows}
of the KP hierarchy is then just an operator $L(t,t_1,t_2,\dots)$ that
gives (formal) integral curves of all these vector fields simultaneously.
Note that
the first KP time $t_1$ is naturally identified with translation along
the original variable $t$.  From the compatibility of the 
equations \eqref{KP flows} in $x=t_2$ and $y=t_3$ (i.e. the fact that 
the corresponding vector fields on the space of Lax operators commute)
one derives that $u=u_1$ satisfies the Kadomtsev--Petviashvili
equation
\begin{equation}\label{KP eq}\frac{3}{4}u_{xx}=(u_y-\frac{1}{4}(6uu_t+u_{ttt}))_t.
\end{equation}
As Sato demonstrated (and we recall in the next section), 
the flows of the KP hierarchy (i.e. our vector fields on the space of
Lax operators) may be reformulated in terms of
a natural abelian Lie
algebra action on (the big cell of) an infinite-dimensional
Grassmannian.  As a consequence, 
the full hierarchy \eqref{KP flows} is easier to understand formally than the
original KP equation \eqref{KP eq}
(though, as the unique solvability of the initial-value problem 
\cite{Mulase e/u} indicates, there is a close relationship 
between solutions $u=u(t,x,y)$ of the single
equation \ref{KP eq} and solutions $L=L(t,t_2,t_3,\dots)$ of the full
KP hierarchy).

\subsection{The Sato Grassmannian}\label{Sato}
For simplicity we concentrate here on the case of the usual KP hierarchy---see 
Section \ref{multicomponent} for comments on the extension to
the multicomponent KP hierarchies.

Consider the vector space $\cV=\C((\del\inv))$ of Laurent series in
$\del\inv$, or constant coefficient microdifferential operators.  It
is convenient to identify $\cV$ with $\cE/\cE\cdot t$, the quotient of
$\cE$ by the left ideal generated by $t$ (i.e. the fiber of $\cE$ at
$t=0$ as a module for $\C[[t]]$ acting by right multiplication). The
vector space $\cV$ has a decomposition 
\bd
\cV=\cV_+\oplus\cV_- = \C[\del]\oplus \del\inv\C[[\del\inv]].
\ed
 Let $\GR(\cV)$ denote the Sato Grassmannian: this
is a space that parametrizes all
subspaces $\cW\subset \cV$ whose projections on 
$\cV/\cV_+$ have finite dimensional
kernel and cokernel (equivalently, $\cW$ is transversal to a subspace
$\cW'$ which is commensurable with $\cV_-$---one thinks of $\cW$ as being
``of the same size as $\cV_+$'').\footnote{In terms of the natural
topology of $\cV$ giving it the structure of {\em Tate vector space}
(\cite{Hecke}), $\GR(\cV)$ parametrizes the $d$-lattices in
$\cV$.} The Sato Grassmannian has a natural structure of a scheme of
infinite type (see e.g. \cite{AMP}), and has connected components
labelled by an integer, the index of a subspace $\cW\subset\cV$ with
respect to $\cV_-$.  A special role is played by an open subset
$\GR^{\circ}(\cV)$ 
of the Sato Grassmannian known as the {\em big cell}: 
this parametrizes subspaces 
$\cW$ that are transversal to $\cV_-$ itself, i.e.
\bd
\GR^{\circ}(\cV)=\big\{\cW\in\GR(\cV): \cW\oplus \C[[\del\inv]] 
=\C((\del\inv))\big\}.
\ed

The algebra $\cE$ acts on $\cV$ by left multiplication. We distinguish
two pieces of the resulting symmetries of $\GR(\cV)$. First, the
constant coefficient microdifferential operators $\Gamma = \C((\del\inv))$ 
act by continuous
endomorphisms of $\cV$, hence give rise to an infinite family of
commuting vector fields on $\GR(\cV)$. Second, consider the {\em
Volterra group}
$$\Emx=\{W=1+w_1\del\inv+w_2\del^{-2}+\cdots\}\subset \cE$$ (where the
$w_i$ are formal power series in $t$).\footnote{Note that 
the Volterra group is a pro-unipotent
algebraic group with Lie algebra $\cE_-\subset \cE$ consisting of
purely negative microdifferential operators.} The group $\cE_-^\times$ acts by
continuous automorphisms on $\cV$, hence algebraically on
$\GR(\cV)$; it also acts on the space $\cL$ of Lax operators via conjugation.
 We will be interested in the action of the commutative
subgroup 
\bd
\Gamx=\cE_-^\times\cap\Gamma = 1+\del\inv\C[[\del\inv]]
\ed
 of constant coefficient
operators (which commutes with the Lie algebra action of $\Gamma$).

\begin{thm}[Sato]
\mbox{}
\begin{enumerate}
\item The action of $\cE_-^{\times}$ on the  big cell 
$\GR^{\circ}(\cV)$ of the
index zero Grassmannian is simply transitive, i.e. every $\cW$ has
 the form $\cW = W\cdot \cV_+$ for a unique
$W\in\cE_-^{\times}$.
\item The action of $\cE_-^{\times}$ on the space $\cL$ of Lax operators
$L=\del+\cdots$ by conjugation is transitive with stabilizer
$\Gamx= 1+\del\inv\C[[\del\inv]]$, i.e. every $L$ is written in the form
$L = W\del W\inv$ with $W\in \cE_-^{\times}$ unique up to constant
coefficient operators.
\item The resulting isomorphism 
\bd
\GR^{\circ}(\cV)/\Gamx \xrightarrow{\sim}\cL
\ed
taking
 $W\cdot\cV_+\mapsto W\del W\inv$
 identifies the infinitesimal
 action of $\del^n\in\C[\del]\subset\Gamma$ on
$\GR^{\circ}/\Gamx$ with the $n$th KP flow $\dfrac{\partial}{\partial
t_n}$ on Lax operators.
\end{enumerate}
\end{thm}
\noindent
The operator $W\in\cE_-^\times$ associated to a subspace or Lax
operator is known as the associated {\em wave operator}.

  Sato's discovery
of the matching of Lax operators with subspaces passes through a
useful intermediary, namely a $\D$-module model of the big cell of
$\GR(\cV_+)$: the big cell is identified with the space of right
$\D$-submodules $M\subset\cE$ that satisfy the transversality property
$\cE=M\oplus\cE_{-}$.  The identification of such 
$\D$-submodules with points of the Sato Grassmannian comes by sending 
$M$ to its fiber at $t=0$,
i.e. the subspace $M/M\cdot t\subset \cV=\cE/\cE\cdot t$ (note that
this subspace no longer carries any $\D$-module structure in general).

A $\D$-submodule $M\subset \cE$ with the transversality property
$\cE =M\oplus \cE_-$
 is automatically a {\em cyclic}
$\D$-module, $M\cong\D$. There is therefore a unique monic
microdifferential operator (``wave operator'') $W\in\Emx$ as above
with the property that $W\cdot M=\D$. 
The Lax operator associated to
$M$ is then $W\del W\inv$.
 One of our
aims in Section \ref{microopers} will be to extend this $\D$-module
description of the big cell to the full Grassmannian.

\section{$\D$-Bundles, Micro-Opers and Noncommutative Geometry}\label{D-bundles}

\subsection{Noncommutative Geometry of $\D$-Modules}\label{quantized}
In this section we describe the noncommutative geometry approach to
$\D$-modules on a curve that motivates most of our constructions;
see also Section \ref{stringy} for
relations with noncommutative gauge theory.
Fix a curve $X$ over an
algebraically closed field of characteristic zero, and let $\D_X$
denote the sheaf of differential operators on $X$.  We assume for simplicity that
$X$ is smooth, although a similar
discussion will apply in the case of nodal or cuspidal curves if
$\D_X$ is replaced by the appropriate ``log'' version; we refer to
\cite{solitons} for more details.

To an algebraic variety one may associate an abelian category, namely
its category of (quasi)coherent sheaves, that encodes the fundamental
geometry of the variety.
In noncommutative algebraic geometry one takes the abelian category\footnote{In fact, 
it may be better to 
replace the category of quasicoherent sheaves on a variety
 by its derived category as a differential graded category,
\cite{Drin}, and take DG categories as the starting point.}
itself as the starting point---see \cite{Staff-ICM, Staff-vdB} for an introduction to
this point of view.  For example, a fundamental construction
of noncommutative algebraic geometry associates, to a sheaf of noncommutative algebras
$\cA$ on a variety, the category of $\cA$-modules.  

Starting with the sheaf $\D_X$ of
differential operators on a smooth curve this construction defines
a noncommutative variety, the category of (quasicoherent) $\D_X$-modules. 
The algebra $\D_X$ is a deformation of the commutative algebra $\Sym(T_X)$,
the sheaf of symbols of differential operators. Since the latter is the
pushforward to $X$ of the algebra of functions on the cotangent bundle
$T^*X$ of $X$, one thinks of $\D_X$
 as the algebra of functions
on a ``quantization of the cotangent bundle of $X$'', a noncommutative
algebraic surface that we will denote by $\spec(\D)=T^*_\hbar X$. 
It turns out that the intuition provided by thinking of 
the category of $\D_X$-modules in this way is an excellent guide to
many interesting questions and useful constructions concerning 
$\D$-modules.

\subsubsection{Completion.}
In studying moduli problems for sheaves, one prefers to work with a
projective (or proper) variety. In particular, moduli problems of
torsion-free sheaves on the open surface $T^*X$ or on its deformation
$T^*_\hbar X$ (i.e. of $\D$-bundles) will not be well behaved, even
if we assume (as we do in the rest of the section) that the curve $X$
is projective. Thus it is convenient to
{\em compactify} to a proper variety by adding a divisor $X_{\infty}$
at infinity and then consider framed sheaves, namely sheaves that
are trivialized at infinity.  Such framing conditions can be
considered as ``asymptotic decay conditions at infinity.''  Indeed,
framing conditions arise
naturally in Yang-Mills theory when one studies connections on a
noncompact 4-manifold that have suitable decay conditions at infinity on
their curvature: such connections often admit extensions to
connections on the complex projective surface that automatically are
trivialized at infinity.

In the case of the vector bundle $T^*X$ there is a standard completion
to a projective bundle $T^*X\subset\ol{T^*X}=\Pp(T^*X\oplus
\Oo)\supset X_\infty$, by adding the curve $X_\infty=\Pp(T^*X)\cong X$
itself at infinity. Thus we consider the projective bundle with
homogeneous coordinate rings $S^{\bullet} = \Sym(T_X\oplus\theo_X) = \oplus_k
(\operatorname{Sym}^{\leq k} T_X)\cdot t^k$.  The category of 
coherent sheaves on $\ol{T^*X} = \on{{\underline Proj}}(S^{\bullet})$
is equivalent to the quotient of the category of finitely generated graded
$S^{\bullet}$-modules modulo its subcategory of bounded modules.
Furthermore, restriction of a coherent sheaf $\widetilde{M}$
on $\ol{T^*X}$ to the curve $X_{\infty} = \on{{\underline Proj}}(\Sym T_X )$ at infinity is
given by taking the graded $S^{\bullet}$-module $M$ to the graded 
$\Sym(T_X)$-module $M/tM$.

The natural quantization of the algebra $S^{\bullet}$, and hence of
the projective completion $\ol{T^*X}$ of $T^*X$, is given by the Rees
ring ${\mathcal R}(\D_X)=\oplus_k \D_X^k\cdot t^k$ of $\D_X$, where $\D_X^k$
denotes differential operators of order at most $k$. A graded module 
over the Rees ring gives a filtered module
over the filtered ring
$\D_X$ and conversely. Since we are interested in such modules modulo
bounded modules, we should consider $\D$-modules which are {\em
eventually} filtered, in other words only sufficiently high filtered
pieces are defined (more precisely, morphisms are only required to
respect the filtration eventually). This leads us to consider the
(derived) category of eventually filtered $\D$-modules as the
completion $\ol{T^*_{\hbar}X}$ of the noncommutative variety
$T^*_\hbar X$. The associated graded ring of $\D$ is just the commutative ring
$\gr(\D)=\Sym T_X$. It follows that our noncommutative
variety is described by adjoining to $T^*_\hbar X$ the {\em
commutative} curve $X_\infty=\proj(\gr(\D))\cong X$. The restriction
of an eventually filtered $\D$-module $M$ to $X_\infty$ is determined by
its asymptotics, i.e. by the coherent sheaf $\gr_N M$ for $N\gg 0$.

\subsubsection{}
The noncommutative variety $\ol{T^*_\hbar X}$ has the same $K$-group
and deRham cohomology (i.e. cyclic homology) as the variety
$\ol{T^*X}$ of which it is a deformation. (In fact Laumon \cite{Lau
filtered} shows that $T^*_\hbar X$, i.e. the derived category of
$\D$-modules, has the same $K$-group as well). Therefore filtered
$\D$-modules have the same numerical invariants as coherent sheaves
on the completed cotangent bundle. In particular it makes sense to
speak of the Chern classes of a filtered $\D$-bundle, and to try to
deform moduli spaces of sheaves with fixed numerical invariants (such
as Hilbert schemes of points) from $\ol{T^*X}$ to $\ol{T^*_\hbar X}$.

It is interesting to note that the ``surface'' $T^*_\hbar X$ exhibits
both one and two-dimensional features, accounting for some of the
peculiarities of noncommutative instantons.  In particular, it has
(cohomological) dimension one, and every torsion-free sheaf on it
($\D_X$-module) is projective, and so may be considered a vector
bundle. However these ``bundles'' are similar to torsion-free sheaves
in the commutative limit, and their moduli provide (flat) deformations
of the moduli spaces of the latter. In particular they carry a
``second Chern class'' (or instanton number) $c_2(M)$, which may be
defined cohomologically, algebraically (using the filtration), or
geometrically (using the associated graded sheaf on $T^*X$).

\subsubsection{Microlocalization.}
In order to interpret the microdifferential Lax
operators of the KP hierarchy geometrically, we will use the interpretation of
microlocalization in terms of the noncommutative space $\ol{T^*_\hbar
X}$---for more on microlocalization, see for example the books
\cite{Kashi,Shap} and the papers \cite{Spring,EO,AvdBvO}.
Recall that the sheaf $\D$ embeds in a
 sheaf of
algebras $\cE$, the sheaf of {\em microdifferential operators}, defined as
follows. First we adjoin to $\D$ (in local coordinates) the formal
inverse $\del\inv$ of a nonvanishing vector field $\del$ satisfying the Leibniz
rule \eqref{Leibniz}, and then we complete with respect to
powers of $\del\inv$, so that in local coordinates $\cE$ is given by
(noncommutative) Laurent series in $\del\inv$ over $\Oo_X$.  The sheaf
$\cE$ is $\Z$-filtered by order in $\del$, extending the filtration
on $\D$, and the completed associated graded algebra is the commutative algebra
given in local coordinates by $\Oo_X((\xi\inv))$, where $\xi$ is the
symbol of $\del$.
 The subsheaf of microdifferential operators of
degree at most zero forms a subalgebra $\cE_-\subset\cE$ which is
complete with respect to the natural topology.

The geometric interpretation of $\cE$ is clear once we note
that the completed associated graded algebra of $\cE$ is naturally identified
with Laurent series along the section at infinity $X_\infty\subset
\ol{T^* X}$ of the compactified cotangent bundle, 
while the associated graded of the subalgebra $\cE_-$ consists of
Taylor series along the same section.  Informally, while $\del$ plays
the role of coordinate ``along the fibers'' on $T^*_\hbar X$, its
inverse $\del\inv$ plays the role of coordinate near the section
$X_\infty$ at infinity, so we may consider $\cE_-$ as functions in the
formal neighborhood of $X_\infty\subset \ol{T^*_\hbar X}$ and $\cE$ as
Laurent series along this section.\footnote{See \cite{PRo} for a discussion
in the context of Kapranov's formal NC schemes.}  In
particular, for a (right) $\D$-module $M$ (sheaf on $T^*_\hbar X$) we
have an associated $\cE$-module $M_\cE=M\ot_\cD \cE$, which is the
restriction of $M$ to the punctured formal neighborhood of $X_\infty$,
while for a filtered $\D$-module (sheaf on $\ol{T^*_\hbar X}$) we
have an associated $\cE_-$-module $M_-=(M_\cE)_{\leq 0}$, the piece
of nonpositive filtration with respect to the induced filtration of
$M_\cE$, which is the restriction of $M$ to the formal neighborhood of
$X_\infty$.

\begin{remark}[Higher-Dimensional Microlocalization]
For a smooth $n$-dimensional variety $X$, microdifferential operators
$\cE$ form a sheaf on the projectivized cotangent bundle
$X_\infty:={\mathbb P}(T^*X)$, which is isomorphic to $X$ for
$n=1$. Here again $X_\infty$ appears as the divisor at infinity in the
projective completion of $T^*X$ or of its noncommutative deformation
$T^*_\hbar X$ on which $\D$-modules live, and microdifferential
operators are again formal Laurent series along this divisor in the
quantized cotangent bundle.
\end{remark}

\subsection{$\D$-Bundles.}

In this section we introduce $\D$-bundles on a
curve.

\begin{defn} \cite{chiral}
A $\D$-bundle $M$ on $X$ is a locally projective, coherent right
$\D_X$-module.
\end{defn}
\noindent
On a smooth (or more generally cuspidal) curve $X$, any torsion-free
$\D_X$-module is locally projective. However a general rank $1$ $\D$-bundle
$M$ on a curve $X$ is not locally free, but only {\em generically} free:
 away from finitely many points of $X$, $M$ is isomorphic to
$\D$.\footnote{Since $\D$ has a skew field of fractions it follows that
locally projective $\D$-modules have well-defined ranks.}

\begin{example}[Right ideals in $\D_X$]  Every right ideal 
in $\D_X$ is torsion-free, hence a $\D$-bundle of rank $1$. However,
ideals in $\D_X$ are typically not locally free. Consider for example
$X=\aline$, so that $\D_X$ is the Weyl algebra $\C\langle
z,\del\rangle/\{\del z-z\del=1\}$.  The right ideal generated by
$z^2$ and $1-z\del$ is not locally free near $z=0$.  However, a
right ideal of $\D_X$ is generically locally free over $\D_X$, and in
fact is equivalent (under rescalings by $\D_X$) to a right ideal in
$\D_X$ which agrees with $\D_X$ generically. Following \cite{CH ideals}
 (see also
\cite{chiral,BW ideals,cusps}), such an ideal
is determined by the (finite) collection of points $x$ at which it
differs from $\D_X$ and choices of subspaces of $\widehat{\theo}_x$ at
those points $x$; for example, the right ideal above corresponds to
the subspace $\boldc + (z^2)$ in $\boldc[[z]]$. (These collections of
subspaces form the {\em ad\`elic Grassmannian}, see Section \ref{adelic}.)
\end{example}

If $M$ is a rank $1$ $\D$-bundle on $X$, we refer to the
finite subset $S$ of $X$ consisting of points $s\in S$ such that $M$ is 
not isomorphic to $\D$ in any neighborhood of $s$, as the set of {\em cusps} of
$M$. The terminology is motivated by the description of $\D$-modules
developed in \cite{CH cusps}, and amplified and generalized in
\cite{cusps}, in terms of coherent sheaves on singular varieties.

In order to construct moduli spaces for $\D$-bundles on $X$, it is
important to ``compactify'' the surface $T^*_\hbar X$ to the
noncommutative $\pline$-bundle $\ol{T^*_\hbar X}$ as explained in
Section \ref{quantized}, in other words to consider filtered
$\D$-bundles (we require the graded components to be vector bundles
on $X$). Moreover, we would like to fix the geometry of the
$\D$-bundle along the (commutative) curve $X_{\infty}$ at infinity in
$\ol{T^*_\hbar X}$, in other words to consider framed $\D$-bundles:

\begin{defn}\label{framing def}
\mbox{}
\begin{enumerate}
\item For a filtered $\D$-module $M$, we denote by $M|_{X_{\infty}}$
the coherent sheaf on $X_\infty=\bproj(\gr(\D))$ associated to the
graded $\gr(\D)$-module $\gr(M)$.
\item Let $V$ denote a vector bundle on $X$. A $V$-framed
$\D$-bundle is a filtered $\D$-bundle $M$ equipped with an
isomorphism $M|_{X_\infty}\to V$.  
\end{enumerate}
\end{defn}

We let $\Mf_{n}(\ol{T^*_\hbar X}, V)$ denote the moduli stack of
$V$-framed $\D$-bundles of second Chern class $n$.  The precise
definition of the moduli stack $\Mf_{n}(\ol{T^*_\hbar X}, V)$ appears
in \cite{solitons}; for a projective curve $X$, $\Mf_{n}(\ol{T^*_\hbar
X}, V)$ is an algebraic stack, and in the case of an elliptic curve
will be described explicitly using the Fourier--Mukai transform in
Section \ref{Fourier}.

\subsection{Micro-Opers and Lax Operators}\label{microopers}

In this section we introduce enhanced versions of $\D$-bundles, the
micro-opers, which give a geometric form to the Lax operators of the
KP hierarchy. On a (commutative) variety, functions act by
$\Oo$-module endomorphisms of any $\Oo$-module. Analogously, we will
need a slight enhancement of the structure of framed $\D$-bundle,
giving the action of a commutative subalgebra of ``functions'' as
endomorphisms. Throughout this section, $(X,\del)$ will denote a curve
with a fixed nowhere-vanishing vector field. Thus $X$ could be any
projective curve of genus at most one, while in higher genus $X$ must
be affine.

\begin{defn}
A micro-oper structure on an $\Oo$-framed $\D$-bundle $M$ is the
data of an $\cE$-module endomorphism $\del_M$ of $M_\cE=M\ot_\D\cE$,
with principal symbol $\del$ with respect to the induced filtration of $M_\cE$.
\end{defn}
 
\begin{remark} See Section \ref{multicomponent}
and \cite{solitons} for discussion of the higher-rank version of
multi-opers, which will give (matrix) Lax operators for multicomponent KP
hierarchies.
\end{remark}

The endomorphism $\del_M$ commutes with the right action of $\cE$, and
is required to perturb the degree of the filtration by one,
$\del_M:(M_{\cE})_i\to(M_{\cE})_{i+1}$, inducing an isomorphism
$\on{gr}(\del_M):\on{gr}_n M_\cE\to \on{gr}_{n+1} M_\cE$ on the graded
pieces compatible with the framing. Thus a micro-oper structure
consists of a lifting of the vector field $\del$ to $M_\cE$, which we
may think of as (part of) a connection,
\footnote{To be more precise, we may think of the action of $\del$ as a 
part of an additional {\em left} action of $\D$ on $M_\cE$.}
that satisfies a strict form of
Griffiths transversality. 

We see that micro-opers are closely
analogous to the {\em opers} introduced by Beilinson and Drinfeld
\cite{Hecke} following Drinfeld and Sokolov \cite{DS}. A $GL_n$-oper
on $X$ is a rank $n$ bundle $V$ with a full flag $\{V_i\}$ and a
connection, which satisfies strict Griffiths transversality:
$\nabla_\del:V_i\to V_{i+1}$, and induces an isomorphism
$V_i/V_{i-1}\to V_{i+1}/V_i$. It is easy to see that $GL_n$-opers are
naturally identified with monic $n$th order differential operators
acting between line bundles on the curve $X$. Thus $GL_n$ opers give a
coordinate-free form of the Lax operators of the $n$th KdV hierarchy,
which is amenable to generalizations to arbitrary reductive groups.
(See Section \ref{affine opers} for more on the parallels of
micro-opers with opers.)

The importance of micro-opers comes from the following theorem,
extending Sato's description of the KP flows on Lax operators on the
disc to micro-opers on a curve $(X,\del)$. Let $\Gamx$ denote the
multiplicative group $1+\del\inv\C[[\del\inv]]\subset\cE(X)$ and
$\Gap=\C[\del]$ considered as an abelian Lie algebra.  For $x\in X$,
we let $\cE_x$ denote the fiber of $\cE$ at $x$ with respect to the
right $\Oo$-module structure, i.e. $\cE_x=\cE/\cE\cdot{\mathfrak
m}_x$.\footnote{This vector space is isomorphic to the vector space 
$\cV =\C((\del\inv))$ from
Section \ref{Sato} using any local coordinate at $x$, in particular
the coordinate coming from the vector field $\del$ near $x$.}

\begin{thm}\label{microopers and GR}\cite{solitons}
 \begin{enumerate}
\item The quotient $\GR(\cE_x)/\Gamx$ of the Sato Grassmannian is in
bijection with micro-opers on the formal disc at $x$.
\item For any $x\in X$, we have a canonical embedding of the space
$\mOp(X)$ of micro-opers on $X$ into the quotient of the Sato
Grassmannian $\GR(\cE_x)/\Gamx$, preserved by the KP flows (given by
the action of $\Gap$ on $\GR(\cE_x)/\Gamx$).
\item The $\D$-bundle $M$ underlying a micro-oper is locally free
near $x$ if and only if the image of $M$ in $\GR(\cE_x)/\Gamx$ is in the
image of the big cell.
\end{enumerate}
\end{thm}

It follows that a micro-oper on $X$ defines (and is determined by) a
canonical KP Lax operator on the dense open subset $U\subset X$ on
which the associated point of the Grassmannian is in the big cell,
i.e. where the $\D$-bundle is locally free. In this sense
the structure of micro-oper gives precise meaning to the
singularities of a microdifferential operator. 
 Thus the space $\mOp(X)$
provides the natural completed phase space for solutions of the KP
hierarchy whose dependence on the first KP time is meromorphic on
$X$. 

\subsubsection{Trivialized Micro-Opers}
The proof of Theorem \ref{microopers and GR} breaks down into two
parts: relating the full Sato Grassmannian with a parameter space of
$\D$-modules, the {\em trivialized micro--opers}, and reducing a
micro-oper to an ``abelian gauge''.

Let $U\subset X$ denote an open subset.
\begin{defn} A {\em trivialization} of a micro-oper $M$ 
on $U$ is the data of an isomorphism $\wt{\eta}:M\ot_\cD \cE\to \cE$
on $X$, or equivalently of a full embedding $\eta:M\hookrightarrow
\cE$ (i.e. an embedding inducing such an isomorphism), compatible with
the framing on $M$.
\end{defn}

Note that a micro-oper $M$ determines an $\cE_-$-module $M_-$, the
zeroth filtered piece of the filtered $\cE$-module $M_\cE=M\ot_\D
\cE$. A trivialization of a micro-oper can equivalently be
described as an isomorphism of $\cE_-$-modules $M_-\to \cE_-$.

From the noncommutative geometry point of view, a trivialization of
a micro-oper $M$ is a {\em formal framing}, i.e. an extension of the
trivialization (framing) of $M$ on the section $X_\infty$ at infinity
to its formal neighborhood. Such extensions form a torsor for the {\em
pro-unipotent} Volterra group $\Emx$ on $X$, and hence exist
on any affine open $U$.

Let $x\in X$. To a trivialized micro-oper $M\hookrightarrow \cE$ on
$X$ we assign its fiber 
$M_x=M/(M\cdot {\mathfrak m}_x)\to \cE_x = \cE/\cE\cdot {\mathfrak m}_x$. 
One checks that
$M_x$ is a $d$-lattice in the vector space $\cE_x$, i.e. a point
of the Sato Grassmannian $\GR(\cE_x)$. This defines a map from
trivialized micro-opers to the Sato Grassmannian. On the formal
disc, we show that this map is an isomorphism of moduli functors
between trivialized micro-opers and the full Sato Grassmannian.
When the $\D$-module $M\hookrightarrow \cE$ is cyclic, we recover
Sato's description of the big cell $\GR^\circ$. This provides the
analog of Theorem \ref{microopers and GR} for trivialized
micro-opers.

It follows from Sato's identification of the big cell with jets of
microdifferential operators that a trivialized micro-oper $(M,\eta)$
on $U$ defines a canonical element of the Volterra group, the {\em
  wave operator} $W_M\in\cE_-^{\times}(U')$, on the open subset
$U'\subset U$ where $M$ is locally free.  We may use
the vector field $\del$ to define a formal coordinate on the formal
disc $\wh{D}_x$ at $x$. This coordinate induces an isomorphism from
the Sato Grassmannian $\GR=\GR(\cV)$ to the Sato Grassmannian at $x$,
$\GR(\cE_x)$, and hence to trivialized micro-opers on $\wh{D}_x$. Thus
any trivialized micro-oper on $X$ defines by restriction to $\wh{D}_x$
a point of $\GR$.

The algebra $\cE(X)$ acts on the left of $\cE$ by endomorphisms of
$\cE$ as a right $\cE$-module
 and hence on the collection of right $\D$-submodules of $\cE$. Thus $\cE(X)$,
considered as a Lie algebra, acts on trivialized micro-opers on $X$.
It is then automatic from Sato's description of the KP flows as the
action of $\C[\del]$ on $\GR$ that the action of the vector field
$\del_n$ restricts to the action of $\del^n\in\C[\del]$ on the image
of $\mOp(X)$ in $\GR$. Note that the first KP time, the action of
$\del$ itself, simply translates infinitesimally along $X$.

\subsubsection{Gauging Micro-Opers}
To complete the proof of Theorem \ref{microopers and GR}, we show that
every micro-oper has a canonical $\Gamx$-orbit of local
coordinatizations on $M$, namely the trivializations in which the
action of $\del_M$ on $M_\cE$ is identified with the left action of
$\del$ on $\cE$. To see this pick any filtered local trivialization
$\wt{\eta}:M_\cE\to\cE$, and note that since $\del_M$ is acting by
right $\cE$-module endomorphisms of $M_\cE$, $\wt{\eta}(\del_M)$ must
act on $\cE$ by left multiplication by an operator of the form
$\del+a_0+a_1\del\inv+\cdots$. Changing the trivialization by the left
action of $\Emx$ on $\cE$, we can conjugate $\wt{\eta}(\del_M)$ to
$\del$, and do so uniquely up to the centralizer of $\del$ in $\Emx$,
namely $\Gamx$. This defines the desired $\Gamx$-orbit of
trivialized micro-oper, which is the inverse to the forgetful map
from trivialized micro-opers. The rest of Theorem \ref{microopers
and GR} follows easily.

\subsection{KdV and Affine Opers} \label{affine opers}
The point of view on micro-opers and the main construction in Theorem
\ref{microopers and GR} are closely parallel to the interpretation in
\cite{BF} of the Drinfeld--Sokolov generalized KdV hierarchies
\cite{DS}. Namely the micro-opers are the KP analogues of the loop
group bundles with connection called {\em affine opers} in \cite{BF}.
We provide a quick overview of the description of the KdV hierarchies
in \cite{BF}, which will not be needed in what follows.

Let $G$ denote a semisimple algebraic group. An affine oper on $X$ is
a $G$-bundle on $X\times\pline$ equipped with a flag along the section
$X\cdot\infty$ and a connection along $X$ that has a prescribed
``transversality'' relation to the flag.  Every affine oper has a
canonical reduction---the {\em Drinfeld--Sokolov gauge}---to an
abelian (Heisenberg) subgroup of the loop group. In this gauge there
is an evident action of the Heisenberg group on the space of affine
opers is evident, and this action defines the KdV hierarchy of
commuting flows.  If we restrict to a disc in $X$ and refine the
reduction to a trivialization of the Heisenberg bundle, we obtain the
(``thick'') loop Grassmannian $\GR_G=LG/LG_-$ associated to $G$, the
moduli of $G$-bundles on $\pline$ trivialized near one point. Thus the
space of affine opers on the disc is identified with the quotient of
the Grassmannian by an abelian group. Moreover affine opers which
satisfy a genericity condition (namely the corresponding $G$-bundle on
$X\times\pline$ is trivial, so we are in the big cell of the loop
Grassmannian) are in natural bijection with $G$-opers \cite{Hecke},
which (for classical groups) are identified with certain differential
operators on $X$. In particular the big cell of the loop Grassmannian
is (an abelian bundle over) the space of KdV Lax operators with formal
power series coefficients.  The reduction to a Heisenberg group is
interpreted as a formal Higgs field, which defines the germ of a
spectral curve. This gives an algebraic approach to the association of
spectral curves (which are formal branched covers of $\pline$) to
these Lax operators, and establishes a natural bijection between the
moduli space of spectral data and that of differential operators.

\begin{remark}[From KP to KdV]
The reduction from KP to KdV is naturally reflected in a reduction
from micro-opers to affine opers: a micro-oper $(M,\del_M)$ on $X$
with the property that the endomorphism $\del_M^n$ of $M_\cE$
preserves the submodule $M\subset M_\cE$, for some positive integer
$n$, gives rise to an affine $GL_n$-oper by using the module structure
of $M$ over polynomials in $\del_M^n$ to define a vector bundle over
$X\times\pline$ and an associated formal spectral curve. Moreover the
big cell conditions on micro-opers and affine opers match up under
this correspondence. In particular, the singularities of KP Lax
operators described by micro-opers reduce to the singularities of KdV
Lax operators described by affine opers.
\end{remark}

\subsection{Micro-Opers and the Ad\`elic Grassmannian}\label{adelic}
Micro-opers also give rise to solutions of the KP hierarchy by a
completely different route, that is closely related to the description
of ideals in the Weyl algebra and more general projective
$\D$-modules developed in \cite{CH ideals,CH cusps,BW
automorphisms,BGK2} (see \cite{cusps} for a different approach, closer
in spirit to the current work, as well as generalizations to higher
dimensions).

It is easy to check (see \cite{solitons}) that a framed $\D$-bundle on $X$
carries a {\em canonical} trivialization (identification with $\D$) on
the open set where it is smooth (i.e. the local data lie in the big
cell).  Note that the existence of this canonical generic
trivialization is to be expected from the point of view of bundles on
a ruled surface. If $S\to X$ is a (commutative) ruled surface with
section $X_\infty\subset S$, then a framed vector bundle on $S$ has a
well-defined generic splitting type. If this generic type is trivial,
then the additional data of a framing determines a unique
trivialization of the bundle away from finitely many jumping points.

The isomorphism classes of $\D$-bundles on a curve $X$ equipped with
a generic trivialization, i.e. embedding as an ideal in rational
differential operators $\D(\C(X))$, are in natural bijection, via the
deRham functor $M\mapsto M\ot_\D\Oo$, with the {\em ad\`elic
Grassmannian} $\Gr^{ad}(X)$ of Wilson \cite{Wilson CM}.
\footnote{More precisely, one obtains the ad\`elic Grassmannian by also
requiring that $M\ot_\D\Oo$ has index $0$ in $\C(X)$ at every point.}
 The latter
parametrizes certain subspaces of the rational functions $\C(X)$ on
$X$, defined by independent conditions at finitely many points of
$X$. Equivalently, $\Gr^{ad}(X)$ parametrizes torsion-free sheaves of
rank one (with generic trivialization) on {\em cuspidal quotients} of
$X$, singular curves having $X$ as their bijective normalization.
(The $\D$-module is obtained from the torsion-free sheaf by a
variant of the induction functor from $\Oo$-modules on $X$ to
$\D$-modules, $M\mapsto M\ot_\Oo\D$.)

In this way we can see that $\Gr^{ad}(X)$ appears as a phase space for
algebro-geometric solutions (\`a la Krichever) of the KP hierarchy,
attached to all cuspidal quotients of $X$. Namely, to a rank one
torsion-free sheaf on a curve, equipped with a trivialization near a
smooth point $\infty$, one assigns its vector space of sections away
from $\infty$, which (using the trivialization and a local coordinate $z\inv$)
define a subspace of $\C((z\inv))$, i.e. a point in the Sato
Grassmannian $\GR(\C((z\inv)))$. 

Thus we have a construction of algebro-geometric solutions of KP from
micro-opers on $X$ (which depends only on the underlying
$\D$-bundle). Note that here $X$ and its cuspidal quotients are
playing the role of the {\em spectral} curve of KP, while in our
construction in Theorem \ref{microopers and GR} the curve $X$
corresponds to the first time of KP, i.e. the (usually formal) curve
where Lax operators live. For a curve of genus greater than one, there
is no intersection between these spaces of ``algebraic'' (Krichever) and
``differential'' (micro-oper) solutions to KP. In the rational, trigonometric and
elliptic cases, it follows from Theorem \ref{the theorem} that all
differential solutions (i.e. micro-opers on cubic curves) are in fact
algebro-geometric solutions, assigned to tangential covers.  In the
rational case, furthermore, there is a symmetry of spectral and
differential variables, namely the Fourier transform, giving rise to
the bispectral involution of Wilson \cite{Wilson bispectral}
identifying the two classes of solutions---see Section
\ref{bispectral}.

\section{Fourier Duality}\label{Fourier}

\subsection{The Fourier--Mukai Transform}
We recall the Fourier--Mukai transform for abelian varieties, in the
special case of an elliptic curve $E$:
\begin{thm} Let $E$ denote an elliptic curve and $\cP$ the Poincar\'e sheaf
on $E\times E$.
\begin{enumerate}
\item (Mukai \cite{Mukai}) The functor $\Ff:M\mapsto Rp_{2*}(p_1^*M\ot\cP)$ defines
an autoequivalence of the bounded derived category of coherent sheaves
on $E$ (and likewise for quasicoherent sheaves).
\item (Laumon \cite{La}, Rothstein \cite{Ro2}) $\Ff$ induces an
equivalence $\Ff:D^b(\D_E)\to D^b(\Enat)$ of bounded derived category
of coherent $\D$-modules on $E$ and of coherent sheaves on $\Enat$.
\footnote{We consider coherent sheaves on $\Enat$ as $\pi_*\Oo_{\Enat}$-modules on $E$.}
\end{enumerate}
\end{thm}

Morally, $\Ff$ describes a coherent sheaf on $E$ as a ``direct
integral'' of degree zero line bundles, which are parametrized by
points of the dual abelian variety, $E$ itself. The extension in the
second part writes a $\D$-module on $E$ as a ``direct integral'' of
flat bundles on $E$, which are parametrized by points of $\Enat$. From
the noncommutative geometry point of view, the first part defines an
automorphism of the ``noncommutative variety'' $D^b(E)$ defined by
$E$, while the second identifies the noncommutative variety
$T^*_{\hbar}E$ with that underlying the commutative variety $\Enat$.

In \cite{solitons} we extend the Fourier--Mukai transform to arbitrary
Weierstrass cubics; see \cite{FriedMorg} for the case of semistable
bundles of degree zero. We also generalize the Fourier--Mukai
transform for $\D$-modules to the singular setting, utilizing results
of \cite{PRo} on Fourier transforms for $D$-algebras. More precisely,
we consider the sheaf $\D^{\log}$ of log-differential operators with
respect to the singularity of the cubic: this is the subsheaf of all
differential operators generated by $\Oo_E$ and the
translation-invariant vector fields coming from the action of the
group $\G$ on $E$. We introduce the surface $\Enat\to E$ as the affine
bundle classifying rank one torsion-free sheaves with a log
connection ($\D^{\log}$-action). 
Reformulating the above theorems in
this setting (with extra care taken along the singularities) we
obtain:

\begin{thm}\label{thm5.2} \cite{solitons} Let $E$ denote a cubic curve and
$\cP$ the Poincar\'e sheaf on $E\times E$.
\begin{enumerate}
\item The functor $\Ff:M\mapsto Rp_{2*}(p_1^*M\ot\cP)$ defines an
autoequivalence of the bounded derived category of coherent sheaves on
$E$ (and likewise for quasicoherent sheaves).
\item $\Ff$ induces an equivalence $\Ff:D^b(\D^{\log}_E)\to
D^b(\Enat)$ of bounded derived category of coherent
$\D^{\log}$-modules on $E$ and of coherent sheaves on $\Enat$.
\end{enumerate}
\end{thm}

The Fourier--Mukai transform is compatible with filtrations on
$\D$-modules and $\Oo_{\Enat}$-modules; as a result, it may be used
(as we will see below) to relate framed $\D$-bundles to coherent sheaves
on $\ol{T^*_{\hbar} E}$.

\begin{remark} It is useful to note that a
degenerate case of the extended Fourier transform gives an
autoequivalence of the derived category of the surface $E\times
\aline$, i.e. of the derived category of modules over the
sheaf of algebras $\Oo_E[s]$. This latter algebra arises as
the common degeneration of
$\D$ and of the coordinate ring $\Oo_{\Enat}$.
\end{remark}

\subsection{Torsion-Free Sheaves and Spectral Sheaves}

In this section we would like to describe the effect of the extended
Fourier--Mukai transform on framed torsion-free sheaves on
$\ol{T^*_\hbar E}$, that is, on framed $\D$-bundles (see Definition
\ref{framing def}). 

\begin{remark} The same
techniques apply to the commutative limit $E\times \pline$, and can be
used to give a new proof of the ADHM classification of framed
torsion-free sheaves on $\pplane$ by quiver data. Namely, we replace
the completion $\pplane$ of $\aplane$ by the completion
$\Gabar\times\pline$, and the Koszul duality and Beilinson
spectral sequence by the Fourier--Mukai transform on $\Gabar$.
\end{remark}

We will consider torsion-free sheaves on $\ol{T^*_\hbar E}$ (and
on $E\times \pline$) trivialized along the section $E_\infty$, and
more generally sheaves framed by a semistable vector bundle $V$ on
$E_\infty$ of degree $0$.  We let $\Mf_{c_2}(S,V)$ denote the moduli
space of these $V$-framed torsion-free sheaves, where
$S=\ol{T^*_\hbar E}$ for $\D$-bundles and $S=E\times\pline$ in the
commutative limit.  We will restrict our attention to framings by
semistable vector bundles $V$ of degree zero whose Fourier transform
$V^{\vee}=\Ff(V)$ is a coherent sheaf on $E$ supported away from the
singular locus. For singular curves $E$, this support condition is
equivalent to the condition (which appeared in Section \ref{Calogero
and Hitchin}) that the pullback of $V$ to the normalization $\pline$
is trivial. 

On the Fourier dual side, we will consider coherent sheaves on
$\Enatbar$ (or in the commutative case $E\times \pline$) whose
restriction to the section $E_\infty$ is identified with a torsion
sheaf $V^{\vee}$ that is supported away from the singular locus. It follows
that such a sheaf has support of dimension at most one.

 A $V^{\vee}$-{\em framed spectral sheaf} on 
$\Enatbar$ or $E\times\pline$ is a coherent sheaf $M$ of pure
dimension one, equipped with an identification
$\phi:M|_{E_{\infty}}\to V^{\vee}$.\footnote{Note that if $E$ is singular
we also add the technical condition that the spectral sheaf is locally 
free over the singular point---see \cite{solitons} for details.}
We similarly define
$V^{\vee}$-framed spectral sheaves on $E\times\pline$. The moduli
space of $V^{\vee}$-framed spectral sheaves on $S^{\vee}=\Enatbar$ or
$S^{\vee}=E\times \pline$ with support of degree $n$ over $E$ is denoted
$\SShv_n(S^{\vee},V^\vee)$.

Note that in the special case $V=\Oo_E$ we have $V^\vee=\Oo_o$ (the
skyscraper at the origin), and so we recover the definition of CM spectral
sheaf on $\Enatbar$ from Definition \ref{CM sheaf}. 
Any $\Oo_o$-framed spectral sheaf is a rank one torsion-free sheaf
on its support $\Sigma\subset\Enatbar$.

The Fourier transform of a $V$-framed torsion-free sheaf is {\it a
priori} a complex of sheaves, whose restriction to $E_\infty$ is
identified with the torsion sheaf $V^{\vee}$. We then prove that the
Fourier transform is itself in fact a sheaf, in cohomological degree
one, and of pure dimension one, and thus defines a $V^\vee$-framed
spectral sheaf. We obtain the following theorem:

\begin{thm} \label{Fourier theorem} \cite{solitons} 
Fix a cubic curve $E$ and a semistable degree zero bundle $V\in
\Bun^{ss}_0(E, n)$.  Then the Fourier--Mukai transform induces an
isomorphism of the moduli spaces $\Mf_{n}(S,V)$ of framed
torsion-free sheaves with $c_2=n$ and $\SShv_{n}(S^{\vee},V^{\vee})$
of framed spectral curves of degree $n$ over $E$, where $S$ is
$\ol{T^*_\hbar E}$ or $E\times\pline$ and $S^{\vee}$ is $\Enatbar$ or
$E\times\pline$, respectively.
\end{thm}

The simplest case of this theorem is the case $V=\Oo_E$,
$V^{\vee}=\Oo_o$. In the commutative case, the component of the moduli
space $\Mf_{n}(E\times\pline,\Oo_E)$ of rank one framed torsion-free
sheaves of sheaves, corresponding to $c_1=0$, is canonically
identified with the Hilbert scheme $(E\times\aline)^{[n]}$ of $n$
points on $E\times\aline$, via the map associating to an ideal of
codimension $n$ in $\Oo_{E\times\pline}$ the underlying torsion-free
sheaf. Thus we obtain a description of the Hilbert scheme of $n$
points on $E\times\aline$ as the space of spectral sheaves framed by
$\Oo_o$.
Since the spectral sheaves on $\Enatbar$ framed by $\Oo_o$ are
precisely the CM spectral sheaves, the noncommutative version of this
leads to the following special case of Theorem \ref{Fourier theorem}:

\begin{corollary}\label{ideals corollary}
The Fourier--Mukai transform induces an isomorphism of schemes
$\Mf_n(\ol{T^*_\hbar E},\Oo)\to\CM_n(E)$.
\end{corollary}

Thus we obtain a description of the {\em completed} phase space of the
Calogero--Moser $n$-particle system as a ``configuration space of $n$
points on the quantized cotangent bundle''. Note that this is a
stronger statement than the natural identification of the
(uncompleted) phase space of {\em distinct} particles with a configuration
space on the cotangent bundle, or more generally the {\em birational}
identification of Hitchin systems with Hilbert schemes of points on
(commutative) cotangent bundles, \`a la Hurtubise \cite{Hu}.

Recall that the Hilbert scheme of $n$ points on $\aplane$ has an
elementary description, as the set of conjugacy classes of pairs of
commuting $n\times n$ matrices $[X,Y]=0$ with a common cyclic vector
$\C[X,Y]\cdot v=\C^n$ (the associated ideal is the kernel of the
projection $\C[X,Y]\to \C^n$). Corollary \ref{ideals corollary} in the
rational (cuspidal) case $E=\Gabar$, combined with the elementary
description of the rational Calogero--Moser phase space, has the
following immediate consequence:

\begin{corollary}\label{CM matrices}
The set of isomorphism classes of finitely generated, rank 1,
torsion-free right modules for the first Weyl algebra $\D_{\aline}$ is
in natural bijective correspondence with the union, over all $n\geq
0$, of the spaces \bd \CM_n=\lbrace (X,Y)\in \gl_n\times \gl_n
\;\big|\; [X,Y] \in \Of\rbrace/GL_n.  \ed
\end{corollary}

This result is immediately implied by combining the two descriptions
of the ad\`elic Grassmannian in \cite{Wilson CM} and \cite{CH ideals},
and was proven in \cite{BW ideals,BGK1,BGK2} using calculations in
noncommutative algebraic geometry. 
Our approach also
gives concrete descriptions of ideals in $\D$ over $\G_m$
and over an elliptic curve, as well as higher rank
versions.

\begin{remark}
What we describe in this section is the Fourier transform for
$\D$-modules on the projective curve $\Gabar$.  There
is another Fourier transform for $\D$-modules on $\aline =
\Ga$, which we discuss in Section \ref{bispectral}, and
which in particular gives a convenient realization of Wilson's
bispectral involution.
\end{remark}

\subsection{The KP/CM Correspondence}\label{final}
As we have mentioned earlier,
an important special class of solutions to the KP equation, first
investigated by Krichever \cite{Kr1,Kr2} and the Chudnovskys
\cite{CC}, consists of functions $u=u(t,x,y)$ which are rational,
trigonometric or elliptic functions of the first KP time $t$,
i.e. extend to rational functions on the additive group $\C$, the
multiplicative group $\Cx$ or an elliptic curve $E$. An analogous
question for the KP hierarchy (studied in \cite{Sh,Wilson CM} in the
rational case) seeks to describe KP Lax operators $L$ which are
rational, trigonometric or elliptic as functions of the first KP time,
namely those $L$ whose orbit under the vector field
$\frac{\partial}{\partial t}$ closes up to (a Zariski open subset of) an
additive, multiplicative or elliptic group.  More generally we would
like to understand the orbits of the KP flows coming from micro-opers
on an arbitrary curve $X$ carrying a nowehere-vanishing vector field
$\del$.

A complete description of the rational, trigonometric and elliptic
solitons of KP follows from the Fourier--Mukai transform, specifically
Corollary \ref{ideals corollary} identifying framed $\D$-bundles with
second Chern class $n$ with the $n$th Calogero--Moser phase space.  We
have identified KP Lax operators with micro-opers, which are framed
$\D$-bundles with an additional endomorphism $\del_M$. However, in
the global setting of a cubic curve there is almost a unique choice of
$\del_M$:

\begin{lemma} A framed $\D$-bundle $M$ 
on a cubic curve carries a canonical micro-oper structure $\del_M$,
and every other choice of $\del_M$ differs by multiplication by
$\Gamx$.
\end{lemma}

Thus we may use Corollary \ref{ideals corollary} to identify
micro-opers, i.e. meromorphic KP Lax operators, with the
Calogero--Moser particle system:

\begin{thm}[\cite{solitons}]\label{the theorem}
Let $E$ denote an arbitrary cubic curve. 
\begin{enumerate}
\item The isomorphism between the moduli
space of
micro-opers on $E$ with $c_2=n$ and the $n$th Calogero--Moser space
$\CM_n(E)$ provided by Corollary \ref{ideals corollary} identifies the
KP flows with the Calogero--Moser flows.
\item This isomorphism identifies the positions of the cusps of a
$\D$-bundle (or poles of a KP Lax operator) and of the corresponding
Calogero--Moser particles.
\end{enumerate}
\end{thm}

The proof of the first part of the
 theorem is immediate from the description of both KP
flows and CM flows as modifications of sheaves along the divisor at
infinity. Indeed, let $\cL$ denote a CM spectral sheaf,
supported on a curve $\Sigma\subset\Enatbar$ with unique, transversal
intersection $o$ with the curve $E_\infty$, and $M$ the micro-oper
corresponding to $\cL$. The CM flows act on $\cL$ through
modifications at the point $o$: restricting $\cL$ to the 
formal punctured neighborhood of $E_\infty$, we obtain a sheaf whose
endomorphism algebra is
isomorphic to the field of
 Laurent series; the polar parts of Laurent series act as infinitesimal
 deformations of
$\cL$.
  This commutative algebra of endomorphisms is identified by the
Fourier--Mukai transform with the algebra of endomorphisms of the
$\cE$-module $M_\cE=M\ot_\D\cE$ (the restriction of $M$ to the
``formal punctured neighborhood'' of the section $E_\infty\subset
\ol{T^*_\hbar E}$). But we have described the KP flows on micro-opers
precisely through the action of these endomorphisms, more specifically
of endomorphisms induced by powers of the endomorphism $\del_M$, so
that the identification of the flows is straightforward.

The statement about locations of the poles is also easy using 
the Fourier--Mukai transform: recall (equation \ref{elliptic
particles}) that a configuration
of $n$ Calogero--Moser particles with distinct positions
is described, in the Hitchin system description, by a spectral
sheaf that pushes forward to the direct sum of line bundles
 $\Oo(q_i-o)$ on $E$.   These line bundles correspond to
the $n$ points $q_i$ under the Fourier transform, and these
points determine the positions at which one creates cusps in
the micro-oper.

\begin{remark} It is instructive to compare the above description of elliptic
solitons with the $\D$-module description of the Krichever
construction due to Nakayashiki and Rothstein \cite{N1,N2,Ro1,Ro2}.
Assume for simplicity that the micro-oper $M$ on $E$ corresponds to a
spectral sheaf which is a line bundle $\cL$ on a smooth spectral curve
$\Sigma\to E$. The dual to the map $\on{Jac}\Sigma\to E$ is a map
$E\to \on{Jac}\Sigma$, which is tangent to $\Sigma\hookrightarrow
\on{Jac}\Sigma$ at the point $x=\Sigma\cap\obar$---in other words
$\Sigma\to E$ is a tangential cover (\cite{TV,DM}). The $\D$-bundle $M$
on $E$ is then the restriction to $E$ of the Krichever $\D$-module
$\on{Krich}(\cL)$ on $\on{Jac}\Sigma$, obtained (as explained in
\cite{Ro2}) as the Fourier--Mukai transform of the Abel--Jacobi
pushforward of $\cL(*x)$ from $\Sigma$ to $\on{Jac}\Sigma$. The cusps
of the $\D$-bundle $M$ correspond to the intersection points of $E$
with the theta divisor of $\on{Jac}\Sigma$.
\end{remark}

\section{Further Topics}\label{further aspects}
In this section we sketch some further applications of the techniques
described in the previous directions.

\subsection{The Multicomponent KP/Spin Calogero--Moser Correspondence}\label{multicomponent}
The description and extension of the KP/CM correspondence outlined in
this paper are worked out in greater generality in the paper
\cite{solitons}. In particular, we extend the correspondence to a
relation between the multicomponent KP hierarchy and the spin
generalizations of the Calogero--Moser system, generalizing and
refining the results of \cite{BBKT,T matrix} in several directions.

The multicomponent KP hierarchy is the matrix generalization of the KP
hierarchy, where we replace the algebra $\cE$ with $\gl_n(\cE)$, the
algebra of $n$ by $n$ matrices over $\cE$, and the Sato Grassmannian
with the $n$-component Grassmanian $\GR(\cV^{\oplus n})$ -- see, for
example, \cite{Mulase, Plaza} for more details. The geometric
interpretation of KP in terms of micro--opers is extended to this
setting by replacing $\D$-line bundles by higher rank $\D$-bundles,
with framing by a general semistable vector bundle of degree
$0$. Higher rank micro--opers are a Higgs refinement of this
structure, namely they carry actions of commutative algebras of matrix
microdifferential operators.  On the Fourier dual side one obtains
general spectral curves in $\Enatbar$, whose geometry is determined by
the framing conditions at infinity. In the case of trivially framed
$\D$--bundles, we obtain precisely the spin generalization of the
Calogero--Moser system, \cite{GH}, and thus obtain an extension of
Wilson's description of pole collisions to a spin CM/multicomponent KP
correspondence, in rational, trigonometric and elliptic settings.  For
other kinds of framings we obtain ``multicolored'' spin
Calogero--Moser systems, describing pole motion of different
reductions of multicomponent KP. 

Among the reductions of multicomponent KP are the 2D Toda lattice
hierarchies. Using a difference analog of the techniques of this paper, we
develop in \cite{solitons2} an extension of the geometric picture described above
to a Toda/Ruijsenaars-Schneider correspondence, extending the results of \cite{KrZab}.

\subsection{Bispectrality}\label{bispectral}
Let us recall the notion of bispectrality of differential operators,
introduced by Duistermaat and Gr\"unbaum \cite{DG}.  Differential
operators $L(t,\del_t)$ and $\Lambda(z,\del_z)$ in two variables $t,z$
are said to be {\em bispectral} if there exists a function $\psi(z,t)$
which is simultaneously a parametric family of eigenfunctions for $L$
and $\Lambda$,
$$L\cdot\psi(z,t)=f(z)\psi(z,t), \hskip.3in \Lambda
\cdot\psi(z,t)=g(t)\psi(z,t)$$ with nonconstant $f,g$.  G. Wilson
\cite{Wilson bispectral} discovered a remarkable symmetry of the
collection of rational solitons (decaying at infinity) of KP, the {\em
bispectral involution}, which expresses the bispectrality of these
solutions. Namely to each point $W$ in Wilson's ad\`elic Grassmannian
(the parameter space of the rational solitons) is
associated a {\em Baker function} $\psi_W(z,t)$, depending on the
spectral parameter $z$ and the KP times $t$, which is a parametric
family of joint eigenfunctions for a commutative ring of differential
operators (this is the ring of functions on the associated spectral
curve). The bispectral involution $W\mapsto b(W)$ is characterized by
the property that it interchanges the spectral parameter with the
first KP time, $\psi_{b(W)}(z,t)=\psi_W(t,z)$ (though it is by no
means clear from this characterization that the desired involution
exists). It follows that the KP solution corresponding to $W$ is
bispectral: it is an eigenfunction for ordinary differential operators
both $z$ and in $t$. In \cite{BW automorphisms}, Berest and Wilson use
the Cannings--Holland bijection between the ad\`elic Grassmannian and
ideals in the Weyl algebra $\D(\aline)$ to give a simpler description
of the bispectral involution and derivation of its main
properties. Namely they identify it with the action of the
antiautomorphism
$$\Ff:\D(\aline)\to \D(\aline), \hskip.2in t\mapsto\del_t,\;
\del_t\mapsto t$$ of the Weyl algebra, i.e. the geometric Fourier
transform composed with the map $t\mapsto -t$.\footnote{See the beautiful
survey articles \cite{Wsurv1,Wsurv2,Wsurv3} for overviews of the work
of Wilson and Berest--Wilson.}

Our description of the rational KP solutions in terms of micro-opers
gives a simple conceptual framework for bispectrality in rational KP. 
The crucial
observation is that micro-opers on $\aline$ define KP solutions in
two independent fashions, as was explained in Section \ref{adelic},
which are interchanged by the geometric Fourier transform. More
precisely, we consider micro-opers on $\pline$ that are locally free
as $\D$-modules near $\infty$---such micro-opers are in bijection with
 points of $\Gr^{ad}$.
 We then define natural embeddings of this space in the
Sato Grassmannians $\GR(\C((\del_t\inv)))$ and $\GR(\C((t\inv)))$.
The first map from micro-opers to the Sato Grassmannian uses Theorem
\ref{microopers and GR}: a micro-oper on $\pline$, when
expanded near $z=0$, gives a KP Lax operator and thus by Sato's
construction a point of the Grassmannian $\GR(\C((\del_t\inv)))$.  The
second map uses a micro-oper to define a torsion-free sheaf on a cuspidal
quotient of $\aline$, which (using natural trivialization data at $\infty$)
defines a subspace of $\C((t\inv))$, i.e. a point in
$\GR(\C((t\inv)))$.  The geometric Fourier transform, which defines an
autoequivalence of the category of $\D$-modules on $\aline$, induces
an involution of this space of micro-opers on $\pline$.  The Fourier
transform also identifies $\C((z\inv))$ and $\C((\del_t\inv))$, their
two Grassmannians, and the collections of KP flows on these
Grassmannians given by the action of $\C[\del_t]$ and $\C[t]$,
respectively. It is easy to see that this identification interchanges
the role of the spectral parameter and of the first KP time, and hence
gives the bispectrality of the rational KP solutions.

\subsection{$\cW_{1+\infty}$ and the Adelic Grassmannian.}\label{factorization}
In \cite{W} we describe vertex algebra structures associated to
$\D$-bundles through the adelic Grassmannian. Namely, for a finite
set $I$, the functor of flat families of $\D$-bundles trivialized
away from $I$ points forms an ind-scheme of ind-finite type over
$X^I$.  The directed system of these spaces over $I$, together with
factorization isomorphisms describing the decomposition of these
spaces with respect to disjoint unions $I\coprod J$, make the adelic
Grassmannian into a factorization ind-scheme
(\cite{chiral}). Moreover the (twisted) delta-functions at the
trivial bundle form a factorization algebra, which we identify with
the $\cW_{1+\infty}$-vertex algebra (the enveloping vertex algebra of
the central extension of the Lie algebra $\D(\K)$). It follows that
the moduli stacks of $\D$-bundles are uniformized by
$\cW_{1+\infty}$-vertex algebra, and that we obtain a localization
for representations of the latter algebra (in other words a definition
of sheaves of twisted conformal blocks). If we replace trivializations
of $\D$-bundles by identifications with some nontrivial $\D$-bundle,
we obtain a new continuous family of chiral algebras, parametrized by
geometry of cusps of $\D$-bundles, all of which localize on
$\D$-moduli spaces. This provides a starting point for the
development of the so-called ``$\cW$-geometry'' from string theory.

\subsection{Hilbert Schemes and Separation of Variables}
The space $\CM_n(E)$ of Calogero--Moser spectral sheaves may be
identified birationally with the Hilbert scheme of $n$ points on
$\Enatbar$, as in \cite{Hu} and in \cite{GNR}, where this description
is interpreted as Sklyanin's separation of variables. A remarkable
feature of the KP/CM correspondence, however, is that on the
$\D$-side the entire space $\CM_n(E)$ is realized biregularly as a
Hilbert scheme of points (moduli of rank one torsion-free sheaves) on
the {\em noncommutative} cotangent bundle of $E$. The description of
positions and ``momenta'' of the cusps of a generic $\D$-bundle $M$
(the latter being coordinates on the one-dimensional Schubert cells
in the adelic Grassmannian), which give the Calogero--Moser particles,
determine a canonical birational identification with a Hilbert scheme
of points on a (commutative) twisted cotangent bundle.  This picture
demonstrates that in order to allow collisions, the proper completion
of the phase space is the noncommutative Hilbert scheme, to which the
flows extend, rather than the commutative one. Also our description of
the noncommutative separation of variables as a Fourier--Mukai
transform establishes the speculation of \cite{GNR} that separation of
variables is a T-duality (see also \cite{KS}).

\subsection{Noncommutative Instantons} 
\label{stringy}

$\D$-bundles on a curve $X$, namely holomorphic bundles on the
quantized cotangent bundle $T^*_\hbar X$ (see Section
\ref{quantized}), are part of the subject of noncommutative gauge
theory (\cite{Nekreview}): one expects a noncommutative version of the
Donaldson--Uhlenbeck--Yau Theorem to identify $\D$-bundles with
noncommutative Yang--Mills instantons. In the rational case, i.e. on
the quantum plane $T^*_{\hbar}\aline$, such a correspondence is
provided {\it a posteriori} by the explicit description of
$\D$-bundles by matrices in Corollary \ref{CM matrices}, and its
higher rank generalizations. Namely, $\D$-bundles on $\aline$ are
classified by the deformed ADHM data, which was shown by Nekrasov and
Schwarz \cite{NekSchwarz} to describe the Yang--Mills instantons on
noncommutative $\R^4$ (see also \cite{Wilson CM,KKO,BraNe,BGK2}).
These noncommutative instantons were proposed as a gauge-theoretic
substitute for the moduli spaces of torsion-free sheaves, resolving
the pointlike instanton singularities of the commutative instanton
moduli spaces on $\pplane$. The moduli of noncommutative instantons
(i.e. $\D$-bundles on $\aline$) are algebraically nontrivial (albeit
diffeomorphic) deformations of the classical instanton moduli spaces.
More generally, bundles on $T^*_\hbar X$ (i.e. $\D$-bundles on $X$)
appear to be a good algebraic model for the systems of D0 branes bound
to a D4 brane by a background $B$-field \cite{NekSchwarz, KS, GNR}.
While the $\D$-bundles are {\em projective}, making possible their
noncommutative instanton interpretation (e.g. \cite{Furuuchi}), they
are not locally free, with singularities (cusps) at special points
which correspond to the position and momenta of the corresponding
Calogero--Moser particles, so that the D0 branes in the D4 brane are
naturally modeled by a many-body system. The noncommutativity of
$T^*_\hbar X$, however, masks these singularities within {\em
  projective} modules.

\bibliographystyle{alpha}
\newcommand{\bl}{/afs/math.lsa.umich.edu/group/fac/nevins/private/math/bibtex/}

\bibliography{\bl baranovsky,\bl bradlow,\bl brosius,\bl
bialynicki-birula, \bl carrell,\bl donagi,\bl friedman,\bl frankel,\bl
ginzburg,\bl gomez,\bl griffiths, \bl grojnowski,\bl hartshorne,\bl
hitchin,\bl hungerford,\bl huybrechts,\bl jardim,\bl kapranov,\bl
kapustin,\bl kirwan, \bl kurke,\bl langton,\bl laumon,\bl lehn,\bl
matsumura,\bl milne, \bl nakajima,\bl nevins,\bl segal,\bl sernesi,\bl
thaddeus,\bl viehweg,\bl wilson}

\end{document}